\documentclass[12pt, oneside, reqno]{amsart}
\usepackage{amsmath,amssymb}
\usepackage{enumitem}
\usepackage{booktabs,array}
\usepackage{microtype}
\usepackage{needspace}
\usepackage{multicol}
\usepackage{etoolbox}
\usepackage{placeins}
\usepackage[textwidth=6in, textheight=9in,marginpar=0.75in]{geometry}
\usepackage[colorlinks=true,linkcolor=blue,citecolor=blue,urlcolor=blue]{hyperref}

\numberwithin{equation}{section}
\newtheorem{theorem}{Theorem}[section]
\newtheorem{conjecture}[theorem]{Conjecture}
\newtheorem{proposition}[theorem]{Proposition}
\newtheorem{lemma}[theorem]{Lemma}
\newtheorem{corollary}[theorem]{Corollary}

\theoremstyle{definition}

\newtheorem{remark}[theorem]{Remark}
\newtheorem{example}[theorem]{Example}

\newcommand{\Z}{\mathbb{Z}}
\newcommand{\Q}{\mathbb{Q}}
\newcommand{\KI}{\textup{\textsf{KnotInfo}}}
\newcommand{\tg}{\widetilde{\gamma}}
\newcommand{\Spinc}{\mathrm{Spin}^c}
\newcommand{\HFhat}{\widehat{\mathit{HF}}}
\makeatletter
\DeclareRobustCommand{\kn}[2]{\in@{a}{#1}\ifin@\textup{#1#2}\else\in@{n}{#1}\ifin@\textup{#1#2}\else
  \ensuremath{#1_{#2}}\fi\fi}
\makeatother
\DeclareMathOperator{\coker}{coker}
\DeclareMathOperator{\rk}{rk}

\makeatletter
\patchcmd{\@settitle}{\bfseries}{\normalfont\scshape}{}
  {\ClassError{amsart}{Could not set the title in small caps}{}}
\patchcmd{\@settitle}{\uppercasenonmath\@title}{}{}
  {\ClassError{amsart}{Could not preserve title capitalization}{}}
\patchcmd{\@setauthors}{\MakeUppercase{\authors}}{{\scshape\authors}}{}
  {\ClassError{amsart}{Could not set the author name in small caps}{}}
\g@addto@macro\@setauthors{\vspace{-1.5ex}
  \begin{center}
    \normalfont\footnotesize
    Center for Geometry and Physics, Institute for Basic Science,\\
    Pohang 37673, South Korea\\
    \href{mailto:seongjinlee@ibs.re.kr}{\nolinkurl{seongjinlee@ibs.re.kr}}
  \end{center}
}
\makeatother

\hypersetup{pdftitle={Cosmetic crossings and the knot complement problem},
  pdfauthor={Seong-Jin Lee}}

\begin{document}

\title[Cosmetic crossings and knot complements]{Cosmetic crossings\\ and the knot complement problem}
\author[Seong-Jin Lee]{Seong-Jin Lee}

\begin{abstract}
We prove the cosmetic crossing conjecture, Problem~1.58 in Kirby's list, and the generalized cosmetic crossing
conjecture for every knot whose branched double cover is an L-space; in particular, both hold for all
alternating and all quasi-alternating knots.  We also prove the oriented
knot complement conjecture, Problem~1.81(D), for all closed $3$-manifolds with finite fundamental group, in particular
for all lens spaces: in such a manifold, every knot whose exterior is irreducible and not a solid torus is determined
by its complement.  Both results rest on a monotonicity theorem for
Heegaard Floer correction terms: if a knot in a rational homology sphere has an isotropic, nonzero homology class, then
the sum of all correction terms strictly increases along its surgeries with coefficients $-1/j$ relative to the
rational longitude, so no nontrivial surgery yields the original oriented manifold.  We also show that every purely
cosmetic band surgery on a knot with an L-space branched double cover is weakly nugatory, which answers a question of
Rogers for these knots.
\end{abstract}


\maketitle

\section{Introduction}\label{sec:intro}

The cosmetic crossing conjecture of X.-S.~Lin \cite[Problem~1.58]{Kir97} asserts that a crossing change on a knot in
$S^3$ that is not nugatory always changes the knot, and the oriented knot complement conjecture
\cite[Problem~1.81(D)]{Kir97}, \cite[Conjecture~6.2]{Gor91} asserts that in a closed oriented $3$-manifold, every knot
whose exterior is irreducible and not a solid torus is determined by its complement (Conjecture~\ref{conj:kc}).  In
this paper we prove the cosmetic crossing conjecture for every knot whose branched double cover is an L-space, in
particular for all alternating and all quasi-alternating knots, together with the corresponding statement for
generalized crossing changes (Theorem~\ref{thm:A} and Corollary~\ref{cor:A}), and we prove the oriented knot complement conjecture for all elliptic
$3$-manifolds, in particular for all lens spaces (Theorem~\ref{thm:E}).  Both results rest on a monotonicity theorem for Heegaard Floer correction terms
(Theorem~\ref{thm:C}): a nontrivial Dehn surgery on a knot whose homology class is isotropic and nonzero never
yields the original oriented manifold.

\subsection{Cosmetic crossing changes}\label{subsec:introccc}
Let $K$ be an oriented knot in $S^3$.  A \emph{crossing disk} for $K$ is an embedded disk $D\subset S^3$ that meets
$K$ transversely in two interior points of opposite sign, and its boundary $\partial D$ is called a \emph{crossing
circle}.  A \emph{crossing change} on $K$ is the operation of performing $\pm1$-surgery on a crossing circle; switching
a crossing in a diagram of $K$ is a crossing change, with $D$ a small disk at the crossing.  The crossing change is
\emph{nugatory} if $\partial D$ bounds an embedded disk in the complement of $K$, in which case it does not change the
isotopy type of $K$.  A crossing change is \emph{cosmetic} if it is not nugatory and the resulting knot is isotopic to
$K$ as an oriented knot.  More generally, for a nonzero integer $q$, a \emph{generalized crossing change} of order $|q|$
on $K$ is the operation of performing $1/q$-surgery on a crossing circle \cite[Section~5.1]{Kal}; it is nugatory or
cosmetic under the same conditions, and for $|q|=1$ it is a crossing change.

\begin{conjecture}[Cosmetic crossing conjecture]\label{conj:ccc}
No knot in $S^3$ admits a cosmetic crossing change.
\end{conjecture}

Conjecture~\ref{conj:ccc} is also known as the nugatory crossing conjecture.  It is known for the unknot
\cite[Theorem~1.4]{ST}, for two-bridge knots \cite{Tor}, for fibered knots \cite{Kal}, for genus one knots with nontrivial
Alexander polynomial \cite{Ito}, extending work of Balm, Friedl, Kalfagianni and Powell \cite{BFKP}, and for special
alternating knots \cite{Bon23}.  Recall that a rational homology sphere $Y$ is an \emph{L-space} if $\HFhat(Y)$ is a
free abelian group of rank $|H_1(Y)|$.  Lidman and Moore proved the conjecture for knots whose branched
double cover $\Sigma(K)$ is an L-space in which each summand of the first homology has square-free order
\cite[Theorem~1.2]{LM}; this class contains the alternating knots of square-free determinant.  Further obstructions
for knots with L-space branched double covers were obtained by Ito \cite{ItoCW} from the Casson--Walker invariant and
by Boninger from the Alexander polynomial: if $\Sigma(K)$ is an L-space and $K$ admits a cosmetic generalized crossing
change, then $\Delta_K$ has a factor $f(t)f(t^{-1})$ with $f\in\Z[t,t^{-1}]$ and $f(-1)\neq\pm1$
\cite[Theorem~1.2]{Bon24}.  Still, the conjecture remained open even for alternating knots: Boninger showed that all but
at most eighteen prime alternating knots with eleven crossings satisfy it \cite[Theorem~4.1]{Bon24}.  The generalized
cosmetic crossing conjecture asserts that no knot in $S^3$ admits a cosmetic generalized crossing change
\cite[Conjecture~1.1]{Rog}; it is known for fibered knots \cite[Theorem~5.2]{Kal}, see also \cite[Theorem~1.2]{Rog}.  We
remove these conditions: both conjectures hold for every knot whose branched double cover is an L-space, and of the
resulting knot we assume only that its branched double cover is homeomorphic to $\Sigma(K)$.  Throughout, $Y\cong Y'$
denotes an orientation-preserving homeomorphism.

\begin{theorem}\label{thm:A}
Let $K$ be a knot whose branched double cover $\Sigma(K)$ is an L-space.  If a generalized crossing change on $K$
yields a knot $K'$ with $\Sigma(K')\cong\Sigma(K)$, for instance a knot isotopic to $K$ as an oriented or as an
unoriented knot, then the generalized crossing change is nugatory.  In particular, $K$ satisfies the cosmetic crossing
conjecture and the generalized cosmetic crossing conjecture.
\end{theorem}

Since the branched double covers of quasi-alternating knots are L-spaces \cite[Proposition~3.3]{OS05}, and alternating
knots are quasi-alternating \cite[Lemma~3.2]{OS05}, we obtain the following.

\begin{corollary}\label{cor:A}
Every quasi-alternating knot, and in particular every alternating knot, satisfies the cosmetic crossing conjecture and
the generalized cosmetic crossing conjecture.
\end{corollary}

The hypothesis of Theorem~\ref{thm:A} covers most of the small knots for which \KI{} \cite{KnotInfo} listed the
conjecture as open in November 2026: of the $214$ such prime knots with at most twelve crossings, $130$ have
branched double covers that are L-spaces, at least with coefficients in $\Z/2$, which suffices
(Section~\ref{subsec:nugatory}), and the branched double covers of the other $84$ are not L-spaces
(Appendix~\ref{app:census}).  For $14$ of the $130$, each summand of $H_1(\Sigma(K))$ has square-free order, so the
theorem of Lidman and Moore applies once $\Sigma(K)$ is known to be an L-space; the other $116$ do not satisfy its
homological hypothesis.

\subsection{Knots determined by their complements}\label{subsec:introkc}
A knot $\kappa$ in a closed oriented $3$-manifold $Y$ is \emph{determined by its complement} if every knot in $Y$
whose complement is homeomorphic to that of $\kappa$ by an orientation-preserving homeomorphism is the image of
$\kappa$ under an orientation-preserving homeomorphism of $Y$.  This holds if no nontrivial Dehn surgery on $\kappa$
yields $Y$ by an orientation-preserving homeomorphism \cite[Section~1]{Gai}.

\begin{conjecture}[Oriented knot complement conjecture]\label{conj:kc}
In a closed oriented $3$-manifold, every knot whose exterior is irreducible and not a solid torus is determined by
its complement.
\end{conjecture}

This is Problem~1.81(D) in Kirby's list \cite{Kir97}, see also \cite[Conjecture~6.2]{Gor91}, in the form stated in
\cite[Conjecture~1.2]{Gai} and \cite[Conjecture~2]{DMR}.  Gordon and Luecke proved it for knots in $S^3$
\cite{GLu89}.  It is also known for null-homologous knots in L-spaces \cite[Theorem~8.2]{Gai}, for null-homotopic
knots in arbitrary closed $3$-manifolds \cite[Theorems~3 and~4]{DMR}, and for non-hyperbolic knots in lens spaces
\cite[Corollary~1.2]{Mat10}; Matignon also states it for non-hyperbolic knots in closed, atoroidal, irreducible Seifert
fibered $3$-manifolds \cite[Theorem~1.3 and Corollary~1.4]{Mat10}.  For all knots in a lens space $L(p,q)$ it was known when $p$ is
square-free \cite[Corollary~8.3]{Gai} and when $p=d^2p'$ with $p'$ square-free and $d\in\{2,3,6\}$ prime to $p'$
\cite[Theorem~1.7(ii)]{ItoCW}, and for all knots in a $3$-manifold with finite fundamental group and square-free first
homology that is obtained by surgery on a knot in $S^3$ \cite[Theorem~1.7(iii)]{ItoCW}.  We prove it for all elliptic
$3$-manifolds, that is, the quotients of $S^3$ by finite groups of isometries acting freely, and for all knots whose
homology class is isotropic and nonzero.  By the geometrization theorem \cite{MT14}, the elliptic $3$-manifolds are
exactly the closed $3$-manifolds with finite fundamental group.

\begin{theorem}\label{thm:E}
\begin{enumerate}[label=(\alph*)]
\item Let $\kappa$ be a knot in a rational homology sphere $Y$ whose homology class is isotropic for the linking form
and nonzero.  Then no nontrivial surgery on $\kappa$ yields $Y$ by an orientation-preserving homeomorphism; in
particular, $\kappa$ is determined by its complement.
\item Conjecture~\ref{conj:kc} holds for all elliptic $3$-manifolds, in particular for all lens spaces: in an elliptic
$3$-manifold, every knot whose exterior is irreducible and not a solid torus is determined by its complement.
\end{enumerate}
\end{theorem}

Here lens spaces are the manifolds $L(p,q)$ with $p\ge1$; this excludes $S^2\times S^1$, and for $S^3$ part~(b) is the
theorem of Gordon and Luecke.  For lens spaces, part~(b) follows from part~(a), the theorem of Gainullin for
null-homologous knots \cite{Gai}, the cyclic surgery theorem \cite{CGLS}, and, for cables of the cores of Heegaard
solid tori, a computation with Seifert invariants \cite{GeL18}.  For the other elliptic manifolds we use, besides part~(a) and Gainullin's theorem,
the finite filling theorem of Boyer and Zhang \cite{BZ01}, the theorem of Culler, Gordon, Luecke and Shalen on essential
tori \cite{CGLS} together with Gordon's description of surgery on cable spaces \cite{Gor83}, the uniqueness of Seifert
fibrations \cite{Hat,Orl72}, the chirality of these manifolds \cite{McC02}, and a linking-form argument for two fillings
at distance three (Lemma~\ref{lem:distance3}).  The orientation is essential.  Bleiler, Hodgson and Weeks found a knot in $L(49,18)$ with a nontrivial
surgery yielding $-L(49,18)$ \cite{BHW}, and Matignon constructed infinitely many knots in lens spaces with nontrivial
surgeries yielding the lens space \cite[Theorem~1.2]{Mat17}, which again reverse orientation (Corollary~\ref{cor:lens}).
By Theorem~\ref{thm:E}(a), every surgery of this kind on a knot whose homology class is isotropic and nonzero reverses
orientation.  The exterior of
the figure-eight sister knot has two fillings that are elliptic but not lens spaces and are homeomorphic only by
orientation-reversing homeomorphisms (Remark~\ref{rem:sister}); for these fillings the orientation is detected by the
linking form (Lemma~\ref{lem:distance3}).

\subsection{Monotonicity of correction terms}\label{subsec:intromono}
Both conjectures lead to the following question: when does a nontrivial Dehn surgery on a knot $\kappa$ in a rational
homology sphere $Y$ yield $Y$ again?  If no such surgery exists, then $\kappa$ is determined by its complement
(Section~\ref{subsec:introkc}).  For a crossing
change on a knot $K$, the Montesinos trick identifies the crossing change with a Dehn surgery on the lift
$\tg\subset\Sigma(K)$ of an arc in $D$ joining the two points of $K\cap D$.  If the crossing change preserves
$\Sigma(K)$, then $\tg$ is isotropic for the linking form, and the surgery has coefficient $\pm1/2$ with respect to the
framing given by the rational longitude.  If $\tg$ is null-homologous and $\Sigma(K)$ is an L-space, then the Dehn
surgery characterization of the unknot of Kronheimer, Mrowka, Ozsv\'ath and Szab\'o \cite{KMOS}, extended by Gainullin
to null-homologous knots in L-spaces \cite{Gai}, together with the equivariant Dehn lemma, shows that the crossing
change is nugatory \cite[Remark~3.4]{LM}.  The remaining case common to both problems is that of a nonzero isotropic
class, and it is settled by the following theorem.  For a rational homology sphere $Y$ we write
\[
S(Y)=\sum_{\mathfrak t\in\Spinc(Y)}d(Y,\mathfrak t)
\]
for the sum of the correction terms of Ozsv\'ath and Szab\'o \cite{OS03}.

\begin{theorem}\label{thm:C}
Let $Y$ be a rational homology sphere, and let $\kappa\subset Y$ be a knot whose homology class is isotropic and has
order $o$.  Let $M$ be the exterior and $\mu$ the meridian of $\kappa$, and let $\lambda$ be the rational longitude of
$M$, oriented so that $\mu\cdot\lambda=1$.  Then for every $j\in\Z$,
\[
S\bigl(M(\mu-(j+1)\lambda)\bigr)-S\bigl(M(\mu-j\lambda)\bigr)\ \ge\ \frac{|H_1(Y)|\,(o^2-1)}{6\,o^2}.
\]
In particular, if $\kappa$ is not null-homologous, then no nontrivial surgery on $\kappa$ yields $Y$ by an
orientation-preserving homeomorphism.
\end{theorem}

Here $\lambda$ is a framing of $\kappa$ because $\kappa$ is isotropic (Lemma~\ref{lem:isotropic}), and
$M(\mu-j\lambda)$ is the result of surgery on $\kappa$ with coefficient $-1/j$ with respect to $\lambda$.  For a
null-homologous knot the bound is zero, and for knots in integral homology spheres Theorem~\ref{thm:C} then recovers
inequalities of Ozsv\'ath and Szab\'o \cite[Corollary~9.14]{OS03}.  The bound is sharp.  A crossing change turns the
stevedore knot $6_1$ into its mirror image; if $\kappa$ is the lift of its crossing arc, then the branched double
covers of the two knots are $M(\mu)$ and $M(\mu\pm2\lambda)$, and equality holds at both steps, from $M(\mu)$ to
$M(\mu\pm\lambda)$ and from $M(\mu\pm\lambda)$ to $M(\mu\pm2\lambda)$.  Equality also holds for the example of Bleiler,
Hodgson and Weeks (Section~\ref{subsec:examples}), and these examples show that Theorem~\ref{thm:C} does not extend to
orientation-reversing homeomorphisms.  The proof applies an inequality of Ozsv\'ath and Szab\'o, which follows from
\cite[Proposition~9.4]{OS03} as in the proof of \cite[Theorem~9.6]{OS03}, to the $2$-handle cobordisms between
consecutive fillings, which are negative definite with intersection form $(-o^2)$, and sums it over a set of
spin$^c$ structures that restricts bijectively to both ends.

Applied to the lift of a crossing arc, Theorem~\ref{thm:C} gives the following statement for every knot.

\begin{theorem}\label{thm:B}
Let $K$ be a knot.  If a generalized crossing change on $K$, for instance a crossing change, yields a knot $K'$ with
$\Sigma(K')\cong\Sigma(K)$, then the lift of the crossing arc to $\Sigma(K)$ is null-homologous.
\end{theorem}

Thus the cosmetic crossing conjecture for an arbitrary knot $K$ reduces to a question about null-homologous knots in
$\Sigma(K)$: it holds for $K$ if every null-homologous knot in $\Sigma(K)$ with a surgery of coefficient $\pm1/2$
yielding $\Sigma(K)$ bounds a disk (Section~\ref{subsec:ccc}).  The same cobordisms were considered by Caudell
\cite{Cau}, whose withdrawn preprint asserted Theorem~\ref{thm:B} for crossing changes that preserve the knot type;
its argument relied on the assertion that, for an integral surgery, the orders of the classes of a knot and of its
dual knot are coprime \cite[Proposition~6]{Cau}.  This fails, for instance, for the integral surgery from
$\Sigma(K_{-1})$ to $\Sigma(K_0)$ in Example~\ref{ex:pretzel}, where $K_{-1}$ is $6_1$ up to mirror image and $K_0$ is
the connected sum of the two trefoils; there both orders are three.

\subsection{Further applications}\label{subsec:introfurther}
The same method applies to band surgeries.  A non-coherent band surgery is a two-strand twist of order one in the
sense of Rogers \cite{Rog}, and a band surgery that yields a knot isotopic to the original one is called purely
cosmetic \cite{IJ}.  A band surgery on a knot that yields a knot is non-coherent, since a coherent band surgery on a
knot yields a link with two components.  Rogers showed that every knot admits purely cosmetic band surgeries that are not nugatory, he
introduced the weaker notion of a weakly nugatory band surgery, and he asked whether some nontrivial knot admits a
purely cosmetic band surgery that is not weakly nugatory \cite[Question~4.7]{Rog}.

\begin{theorem}\label{thm:D}
Let $K$ be a knot whose branched double cover is an L-space.  If a non-coherent band surgery on $K$ yields a knot $K'$
with $\Sigma(K')\cong\Sigma(K)$, then the band surgery is weakly nugatory.  In particular, every purely cosmetic band
surgery on $K$ is weakly nugatory.
\end{theorem}

Since $\Sigma(K')\cong\Sigma(K)$ for mutant knots $K$ and $K'$ \cite{Vir76},
Theorem~\ref{thm:A} implies that two mutant knots with L-space branched double covers that are not isotopic are not
related by a crossing change (Corollary~\ref{cor:mutants}).  This fails for non-coherent band surgeries: a nugatory
band surgery can turn a knot into a mutant that is not isotopic to it (Remark~\ref{rem:bandmutants}).
Finally, together with a theorem of Ito \cite[Corollary~1.11]{ItoCW}, Corollary~\ref{cor:A} gives the following.

\begin{corollary}\label{cor:B}
Every prime knot with at most ten crossings satisfies the cosmetic crossing conjecture.
\end{corollary}

This paper is organized as follows.  Section~\ref{sec:background} collects background on correction terms, rational
longitudes, the Montesinos trick, and the criteria of Lidman and Moore and of Rogers.  Section~\ref{sec:monotone}
contains the proof of Theorem~\ref{thm:C} and the examples.  Section~\ref{sec:applications} contains the
applications: Theorems~\ref{thm:A} and~\ref{thm:B} for crossing changes, Corollary~\ref{cor:B}, and their extension to
links (Section~\ref{subsec:ccc}), generalized crossing changes (Section~\ref{subsec:gcc}), Theorem~\ref{thm:D} and
mutants (Section~\ref{subsec:bands}), and Theorem~\ref{thm:E} (Section~\ref{subsec:complement}).  In Appendix~\ref{app:alternating} we give a second proof of Corollary~\ref{cor:A}
for alternating knots, which replaces Theorem~\ref{thm:C} by Donaldson's diagonalization theorem, and in
Appendix~\ref{app:census} we list the knots with at most twelve crossings for which \KI{} listed the conjecture as open
in November 2026, according to whether their branched double covers are L-spaces.

\section{Background}\label{sec:background}

In this section we fix conventions and collect the facts on which the rest of the paper is based: the behavior of
correction terms under a negative definite $2$-handle cobordism, the rational longitude of a knot exterior, the
Montesinos trick for crossing changes and two-strand twists, and the criteria of Lidman and Moore and of Rogers for such
an operation to be nugatory.

Throughout, manifolds are compact and oriented, and homology groups have integer coefficients unless indicated
otherwise.  For an oriented manifold $Y$ we denote by $-Y$ the same manifold with the opposite orientation, and for
closed $3$-manifolds we write $Y\cong Y'$ if there is an orientation-preserving homeomorphism from $Y$ to $Y'$.
Branched double covers of $S^3$ carry the orientation lifted from $S^3$.

Let $\kappa$ be a knot in a closed $3$-manifold $Y$, let $\nu(\kappa)$ be an open tubular neighborhood of $\kappa$,
and let $M=Y\setminus\nu(\kappa)$ be the exterior, oriented as a submanifold of $Y$.  A \emph{slope} on $\partial M$
is the isotopy class of an essential simple closed curve; we identify slopes with primitive elements of
$H_1(\partial M)$ up to sign.  The \emph{distance} $\Delta(\eta,\xi)$ between two slopes is their minimal geometric
intersection number, that is, the absolute value of the algebraic intersection number of the corresponding classes.
For a slope $\eta$, we denote by $M(\eta)$ the Dehn filling of $M$ along $\eta$, oriented so as to extend the
orientation of $M$, and we call the core of the filling solid torus the \emph{dual knot} of the filling; its
exterior is $M$ and its meridian is $\eta$.  The \emph{meridian} $\mu$ of $\kappa$ is the slope that bounds a disk in
$\nu(\kappa)$, and
a \emph{framing} of $\kappa$ is a slope $\lambda$ with $\Delta(\mu,\lambda)=1$.  Given a framing, we orient $\mu$ and
$\lambda$ so that $\mu\cdot\lambda=1$ on $\partial\nu(\kappa)$, oriented as the boundary of $\nu(\kappa)$, and we say
that the slope $a\mu+b\lambda$ has \emph{coefficient} $a/b$ with respect to $\lambda$.  For a knot in $S^3$ and its
Seifert framing these are the usual surgery coefficients, and $Y_r(\kappa)$ denotes the result of surgery on $\kappa$
with coefficient $r$ when the framing is understood.

\subsection{Correction terms and negative definite cobordisms}\label{subsec:dinv}
Let $Y$ be a rational homology sphere.  The set $\Spinc(Y)$ of spin$^c$ structures on $Y$ is an affine space over
$H^2(Y)\cong H_1(Y)$, so it has $|H_1(Y)|$ elements.  For $\mathfrak t\in\Spinc(Y)$, Ozsv\'ath and Szab\'o defined
Heegaard Floer homology groups $\HFhat(Y,\mathfrak t)$, $HF^\infty(Y,\mathfrak t)$ and $HF^+(Y,\mathfrak t)$; the
latter two are modules over $\Z[U]$, related by a natural $U$-equivariant map
\[
\pi\colon HF^\infty(Y,\mathfrak t)\longrightarrow HF^+(Y,\mathfrak t)
\]
\cite{OS04a}.  These groups carry absolute $\Q$-gradings, which $\pi$ preserves and $U$ lowers by two
\cite[Theorem~7.1]{OS06}, and $HF^\infty(Y,\mathfrak t)\cong\Z[U,U^{-1}]$ \cite[Theorem~10.1]{OS04b}.
The \emph{correction term} $d(Y,\mathfrak t)\in\Q$ is the minimal grading of a non-torsion element in the image of
$\pi$ \cite[Definition~4.1]{OS03}.  By \cite[Proposition~4.2]{OS03},
\begin{equation}\label{eq:dsym}
d(Y,\bar{\mathfrak t})=d(Y,\mathfrak t)\qquad\text{and}\qquad d(-Y,\mathfrak t)=-d(Y,\mathfrak t),
\end{equation}
where $\bar{\mathfrak t}$ denotes the conjugate of $\mathfrak t$, and $\Spinc(-Y)$ is identified with $\Spinc(Y)$.
Since every homeomorphism between $3$-manifolds is isotopic to a diffeomorphism, an orientation-preserving
homeomorphism $f\colon Y\to Y'$ induces a bijection $f_*\colon\Spinc(Y)\to\Spinc(Y')$ with $d(Y',f_*\mathfrak
t)=d(Y,\mathfrak t)$ for all $\mathfrak t$.  We shall work with the sum of all correction terms,
\begin{equation}\label{eq:S}
S(Y)=\sum_{\mathfrak t\in\Spinc(Y)}d(Y,\mathfrak t).
\end{equation}
It depends only on the orientation-preserving homeomorphism type of $Y$, and $S(-Y)=-S(Y)$ by~\eqref{eq:dsym}.

Let $Y$ and $Y'$ be rational homology spheres.  A \emph{cobordism} from $Y$ to $Y'$ is a smooth, compact, connected
$4$-manifold $W$ with $\partial W=-Y\sqcup Y'$.  Since $H_1(\partial W;\Q)=H_2(\partial W;\Q)=0$, the map
$H_2(W;\Q)\to H_2(W,\partial W;\Q)$ is an isomorphism, and the intersection form of $W$ is nondegenerate on
$H_2(W;\Q)$; we denote by $b_2^+(W)$ and $b_2^-(W)$ the numbers of its positive and negative eigenvalues and by
$\sigma(W)$ its signature.  For $\mathfrak s\in\Spinc(W)$, the square $c_1(\mathfrak s)^2\in\Q$ is defined by means
of the isomorphism $H^2(W,\partial W;\Q)\to H^2(W;\Q)$.  A spin$^c$ structure $\mathfrak s$ on $W$ induces
homomorphisms
\[
F^\infty_{W,\mathfrak s}\colon HF^\infty(Y,\mathfrak s|_Y)\to HF^\infty(Y',\mathfrak s|_{Y'}),\qquad
F^+_{W,\mathfrak s}\colon HF^+(Y,\mathfrak s|_Y)\to HF^+(Y',\mathfrak s|_{Y'}),
\]
which satisfy $\pi\circ F^\infty_{W,\mathfrak s}=F^+_{W,\mathfrak s}\circ\pi$ \cite[Theorem~1.1]{OS06} and shift the
absolute grading by
\begin{equation}\label{eq:shift}
\frac{c_1(\mathfrak s)^2-2\chi(W)-3\sigma(W)}{4}
\end{equation}
\cite[Theorem~7.1]{OS06}; see also \cite[Equation~(4)]{OS03}.  We use these maps only through the following
inequality.

\begin{proposition}\label{prop:dineq}
Let $Y$ and $Y'$ be rational homology spheres, and let $W$ be a cobordism from $Y$ to $Y'$ obtained from
$Y\times[0,1]$ by attaching a single $2$-handle along $Y\times\{1\}$.  If $b_2^+(W)=0$, then
\begin{equation}\label{eq:dineq}
d(Y',\mathfrak s|_{Y'})\ \ge\ d(Y,\mathfrak s|_Y)+\frac{c_1(\mathfrak s)^2+1}{4}\qquad\text{for every }\mathfrak
s\in\Spinc(W).
\end{equation}
\end{proposition}

\begin{proof}
By the exact sequence of the pair $(W,Y)$, we have $H_2(W;\Q)\cong H_2(W,Y;\Q)\cong\Q$.  Hence $b_2^-(W)=1$ and
$\sigma(W)=-1$.  Moreover $\chi(W)=\chi(Y\times[0,1])+1=1$, so the grading shift~\eqref{eq:shift} equals
$(c_1(\mathfrak s)^2+1)/4$.  In the notation of \cite[Section~9]{OS03}, $W$ is the cobordism $W(K)$ associated with
the attaching circle $K\subset Y$ of the $2$-handle and its framing.  Since $HF^\infty(Y,\mathfrak
t)\cong\Z[U,U^{-1}]$ for every $\mathfrak t$, the group $HF^\infty(Y)$ is standard in the sense of
\cite[Section~9]{OS03}, and every spin$^c$ structure on $Y$ or $Y'$ is torsion.  Therefore $F^\infty_{W,\mathfrak s}$
is an isomorphism for every $\mathfrak s\in\Spinc(W)$ by \cite[Proposition~9.4]{OS03}.

Fix $\mathfrak s$, put $\mathfrak t=\mathfrak s|_Y$ and $\mathfrak t'=\mathfrak s|_{Y'}$, and let $y$ be a
homogeneous non-torsion element of grading $d(Y',\mathfrak t')$ in the image of $\pi\colon HF^\infty(Y',\mathfrak
t')\to HF^+(Y',\mathfrak t')$.  Since $F^\infty_{W,\mathfrak s}$ is surjective, and since it and $\pi$ are
homogeneous, there is a homogeneous element $z\in HF^\infty(Y,\mathfrak t)$ with $\pi(F^\infty_{W,\mathfrak
s}(z))=y$.  The element $x=\pi(z)$ is homogeneous and lies in the image of $\pi$, and $F^+_{W,\mathfrak s}(x)=y$.
Since $y$ is non-torsion, so is $x$, and hence the grading of $x$ is at least $d(Y,\mathfrak t)$.  Therefore
\[
d(Y',\mathfrak t')=\operatorname{gr}(y)=\operatorname{gr}(x)+\frac{c_1(\mathfrak s)^2+1}{4}\ \ge\ d(Y,\mathfrak
t)+\frac{c_1(\mathfrak s)^2+1}{4}.\qedhere
\]
\end{proof}

The argument is that of the proof of \cite[Theorem~9.6]{OS03}; the same inequality appears in the proof of
\cite[Proposition~9.9]{OS03}.

\subsection{Rational longitudes and isotropic knots}\label{subsec:isotropic}
Let $M$ be a compact, connected $3$-manifold with torus boundary such that $H_1(M;\Q)\cong\Q$, for instance the
exterior of a knot in a rational homology sphere.  The image of $i_*\colon H_1(\partial M;\Q)\to H_1(M;\Q)$ is
one-dimensional, and there is a unique slope $\lambda_M$ on $\partial M$, the \emph{rational longitude}, such that
$i_*(\lambda_M)$ has finite order in $H_1(M)$ \cite[Definition~3.1]{Wat12}, \cite[Section~2.2 and
Definition~2.2]{LM}.  In particular $i_*(\eta)$
has infinite order for every slope $\eta\neq\lambda_M$.  By a lemma of Watson \cite[Lemma~3.2]{Wat12}, see
also \cite[Lemma~2.3]{LM}, there is a positive integer $c_M$ such that
\begin{equation}\label{eq:watson}
|H_1(M(\eta))|=c_M\,\Delta(\eta,\lambda_M)\qquad\text{for every slope }\eta\neq\lambda_M.
\end{equation}

Let $Y$ be a rational homology sphere.  Its \emph{linking form} $\ell_Y\colon H_1(Y)\times H_1(Y)\to\Q/\Z$ is defined
as follows: if $\xi$ and $\zeta$ are disjoint $1$-cycles representing $x$ and $y$, and $C$ is a $2$-chain transverse
to $\zeta$ with $\partial C=n\,\xi$ for some integer $n\neq0$, then $\ell_Y(x,y)=(C\cdot\zeta)/n\bmod1$.  It is a
nondegenerate symmetric bilinear form.  A class $x\in H_1(Y)$ is \emph{isotropic} if $\ell_Y(x,x)=0$, and a knot in
$Y$ is \emph{isotropic} if its homology class is.  Null-homologous knots are isotropic.  If $|H_1(Y)|$ is
square-free, they are the only isotropic knots, since the linking form is the orthogonal sum of its restrictions to the
primary parts of $H_1(Y)$, and each of these is then cyclic of prime order.  The following elementary lemma
characterizes isotropic knots by their rational longitudes.

\begin{lemma}\label{lem:isotropic}
Let $\kappa$ be a knot in a rational homology sphere $Y$, with exterior $M$ and meridian $\mu$.
\begin{enumerate}[label=(\alph*)]
\item $\Delta(\mu,\lambda_M)=1$ if and only if $\kappa$ is isotropic.
\item If $\Delta(\mu,\lambda_M)=1$, then $[\kappa]\in H_1(Y)$ and $i_*(\lambda_M)\in H_1(M)$ have the same order.
\item If $\kappa$ is isotropic, then for every slope $\eta$ with $\Delta(\eta,\lambda_M)=1$, the manifold $M(\eta)$
is a rational homology sphere with $|H_1(M(\eta))|=|H_1(Y)|$, and its dual knot is isotropic, with homology class of
the same order as $[\kappa]$.
\end{enumerate}
\end{lemma}

\begin{proof}
Let $V=Y\setminus\operatorname{int}M$ be the solid torus with core $\kappa$, let $\rho$ be a framing of $\kappa$, and
write $\lambda_M=a\mu+b\rho$ in $H_1(\partial M)$ with $\gcd(a,b)=1$, so that $\Delta(\mu,\lambda_M)=|b|$; here
$b\neq0$, since $H_1(M(\mu))=H_1(Y)$ is finite.  Let $k\ge1$ be the order of $i_*(\lambda_M)$ in $H_1(M)$.
Regarding $\mu$ and $\rho$ as simple closed curves on $\partial M$, the $1$-cycle $ka\,\mu+kb\,\rho$ is
null-homologous in $M$, so it is the boundary of a $2$-chain $F$ in $M$.  In $V$, the curve $\mu$ bounds a
meridian disk $E$, and $\rho$ cobounds an annulus $R$ with $\kappa$; orient them so that $\partial E=\mu$ and $\partial
R=\rho-\kappa$.  Then the $2$-chain $C=F-ka\,E-kb\,R$ satisfies $\partial C=kb\,\kappa$.  Let
$\kappa^\sharp\subset\operatorname{int}V$ be a parallel copy of $\kappa$ that is disjoint from $R$.  It is disjoint from
$F$ and meets $E$ in one point, so $C\cdot\kappa^\sharp=\pm ka$, and
\[
\ell_Y([\kappa],[\kappa])=\frac{C\cdot\kappa^\sharp}{kb}\bmod1=\pm\frac ab\bmod1 .
\]
Since $\gcd(a,b)=1$, this vanishes if and only if $|b|=1$, which proves~(a).

For~(b), the hypothesis says that $\lambda_M$ is a framing of $\kappa$, so it is isotopic to $\kappa$ in $V$, and
$[\kappa]$ is the image of $i_*(\lambda_M)$ under the quotient map $H_1(M)\to H_1(M)/\langle i_*(\mu)\rangle=H_1(Y)$.
Since $i_*(\lambda_M)$ has finite order and $i_*(\mu)$ has infinite order, the subgroups they generate intersect
trivially, and hence $i_*(\lambda_M)$ and its image $[\kappa]$ have the same order.

For~(c), we have $\Delta(\mu,\lambda_M)=1$ by~(a), so~\eqref{eq:watson} gives
$|H_1(M(\eta))|=c_M=|H_1(Y)|$.  The dual knot of $M(\eta)$ has exterior $M$ and meridian $\eta$, so it is isotropic
by~(a), and by~(b), applied to it and to $\kappa$, both homology classes have the order of $i_*(\lambda_M)$.
\end{proof}

Let $L\subset S^3$ be a framed link with $n$ components whose linking matrix $N$ is nonsingular.  We denote by $X_N$
the $4$-manifold obtained from $B^4$ by attaching $2$-handles along $L$; the notation suppresses the dependence on
$L$.  Then $H_2(X_N)\cong\Z^n$, with intersection form given by $N$ in the basis of the $2$-handles, and
$Y=\partial X_N$ is a rational homology sphere with $H_1(Y)\cong\coker N$, generated by the meridians of the
components of $L$.  By the theorem of Lickorish and Wallace, every rational homology sphere arises in this way.  Let
$\kappa\subset Y$ be a knot.  After an isotopy we may assume that $\kappa\subset S^3\setminus L$, and we call the vector
$v\in\Z^n$ of linking numbers of
$\kappa$ with the components of $L$ a \emph{linking vector} of $\kappa$.  Then $[\kappa]=[v]$ in $\coker N$, so the
order of $[\kappa]$ is the least common multiple of the denominators of the entries of $N^{-1}v$.  We describe slopes
on $\partial\nu(\kappa)$ by their coefficients with respect to the Seifert framing of $\kappa$ in $S^3$, and we set
\[
x_0=v^TN^{-1}v\in\Q .
\]

\begin{lemma}\label{lem:presentation}
The rational longitude of the exterior of $\kappa$ in $Y$ has coefficient $x_0$.  In particular, $\kappa$ is
isotropic if and only if $x_0\in\Z$.
\end{lemma}

\begin{proof}
Let $a/b$ be a slope with $\gcd(a,b)=1$ and $b\ge0$.  The group $H_1(Y_{a/b}(\kappa))$ is generated by the meridians
of $L$ and of $\kappa$, with relations given by the rows of the matrix
\[
\begin{pmatrix}N&v\\ b\,v^T&a\end{pmatrix},\qquad\text{whose determinant is }\det N\cdot(a-b\,x_0).
\]
Hence $H_1(Y_{a/b}(\kappa))$ is infinite exactly when $a/b=x_0$.  By~\eqref{eq:watson}, the rational longitude is
the only slope whose filling has infinite first homology, so it has coefficient $x_0$.
Its distance from the meridian, which has coefficient $1/0$, is the denominator of $x_0$, and the second assertion
follows from Lemma~\ref{lem:isotropic}(a).
\end{proof}

\subsection{Crossing changes, two-strand twists and the Montesinos trick}\label{subsec:montesinos}
Let $K$ be a knot in $S^3$, and let $D$ be a crossing disk for $K$.  Every crossing change on $K$ can be realized by
switching a crossing $c$ in a suitable diagram.  Indeed, if $D'\subset D$ is a regular neighborhood of an arc in $D$
joining the two points of $K\cap D$, then $\partial D$ and $\partial D'$ cobound an annulus in $D\setminus K$; hence
surgeries on $\partial D$ and on $\partial D'$ give the same knot, $\partial D$ bounds a disk in the complement of $K$
if and only if $\partial D'$ does, and $D'$ is a crossing disk at a crossing of a suitable diagram of $K$.  Following
\cite{LM}, we denote by $K_+$ and $K_-$ the knots whose diagrams agree away from $c$ and in which $c$ is a positive
and a negative crossing, respectively; one of them is the given knot.  A \emph{crossing arc} for $c$ is an arc
$\gamma\subset D$ joining the two points of $K_+\cap D$.  Let $Y=\Sigma(K_+)$ be the double cover of $S^3$ branched
along $K_+$.  The arc $\gamma$ lifts to a knot $\tg\subset Y$, and we denote its exterior by $M=Y\setminus\nu(\tg)$,
oriented as a submanifold of $Y$.  By the Montesinos trick \cite[Lemma~2.1]{LM}, there are slopes $\alpha$ and
$\beta$ on $\partial M$ such that
\begin{equation}\label{eq:montesinos}
M(\alpha)=\Sigma(K_+),\qquad M(\beta)=\Sigma(K_-),\qquad \Delta(\alpha,\beta)=2 .
\end{equation}
Both fillings contain $M$ as the double cover of the complement of a ball around $c$, branched along the knot, so the
first two identifications respect orientations; the slope $\alpha$ is the meridian of $\tg$.

The following lemma is contained in the proofs of Theorem~2.4 and Proposition~3.3 of \cite{LM}, where the case $m=2$
below is excluded because the determinant of a knot is odd.  We give an argument that applies to links as well.

\begin{lemma}[Lidman--Moore]\label{lem:LM}
If $\det K_+=\det K_-$, then $\Delta(\alpha,\lambda_M)=\Delta(\beta,\lambda_M)=1$.  Moreover, if $\alpha$ and
$\lambda_M$ are oriented so that $\alpha\cdot\lambda_M=1$ on $\partial\nu(\tg)$, then
$\beta=\pm(\alpha+2\varepsilon\lambda_M)$ for some $\varepsilon\in\{1,-1\}$.
\end{lemma}

\begin{proof}
Since $|H_1(\Sigma(K_\pm))|=\det K_\pm$, it follows from~\eqref{eq:watson} that
$\Delta(\alpha,\lambda_M)=\Delta(\beta,\lambda_M)$; we denote this number by $m$.  Choose a slope $\rho$ with
$\rho\cdot\lambda_M=1$ on $\partial\nu(\tg)$, and write $\alpha=m\rho+r\lambda_M$ and $\beta=\pm(m\rho+s\lambda_M)$.
Then $2=\Delta(\alpha,\beta)=m\,|r-s|$, so $m\in\{1,2\}$.  If $m=2$, then $r$ and $s$ are odd, because $\alpha$ and
$\beta$ are primitive, and so $|r-s|\neq1$, a contradiction.  Hence $m=1$ and $|r-s|=2$, and the choice $\rho=\alpha$
gives the second assertion.
\end{proof}

Thus, whenever $\det K_+=\det K_-$, the knot $\tg$ is isotropic by Lemma~\ref{lem:isotropic}(a), and with respect to
its framing $\lambda_M$ the slope $\beta$ has coefficient $\varepsilon/2$.  Reversing the orientation of $M$ reverses
that of $\partial M$, and the coefficient of $\beta$ becomes $-\varepsilon/2$.

\begin{remark}\label{rem:links}
The Montesinos trick~\eqref{eq:montesinos} is a statement about the double cover of a ball around $c$ branched along
the two arcs in it, so it holds for a crossing change between two links $L_+$ and $L_-$ as well.  If $\det L_+=\det
L_-\neq0$, then $\Sigma(L_\pm)$ are rational homology spheres with
$|H_1(\Sigma(L_\pm))|=\det L_\pm$, and the proof of Lemma~\ref{lem:LM} applies verbatim.
\end{remark}

Rogers studied the following generalization of crossing changes \cite{Rog}.  A \emph{twisting circle} for a knot
$K\subset S^3$ is a curve $c\subset S^3\setminus K$ bounding a disk $D$, a \emph{twisting disk}, that meets $K$
transversely in two points of opposite sign, and an arc $\gamma\subset D$ joining these two points is a
\emph{twisting arc}.  We fix a product neighborhood $B=D\times[0,1]$ of $D$ such that each disk $D\times\{t\}$ meets
$K$ in two points of opposite sign.  For a nonzero integer $n$, the \emph{two-strand $n$-twist} along $D$ modifies
$(B,B\cap K)$ by the twisting homeomorphism of \cite[Definition~2.6]{Rog}, which adds $n$ half-twists to the two
strands of $K$ in $B$; its \emph{order} is $|n|$.  For even $n$ it is the generalized crossing change given by
$-2/n$-surgery on $c$ \cite[Remark~2.7]{Rog}, so it depends only on $c$.  For odd $n$ it depends on $D$, and not
only on $c$ (Remark~\ref{rem:bandmutants}).  For $n=\pm1$
it is a non-coherent band surgery \cite[Section~1]{Rog}, and conversely, every non-coherent band surgery on $K$ can be
realized as a two-strand twist of order one \cite[Remark~4.2]{Rog}.  The twist is \emph{nugatory} if $c$ bounds a disk
in $S^3\setminus K$.  For even $n$ the resulting knot is then isotopic to $K$, since the twist is a surgery on $c$.  For
odd $n$ this need not hold: a nugatory twist of order one can turn $J\#J$ into $J\#rJ$, where $rJ$ is $J$ with its
orientation reversed, and these knots differ when $J$ is not invertible (Remark~\ref{rem:bandmutants}).

As before, a twisting arc $\gamma$ lifts to a knot $\tg\subset\Sigma(K)$.  Let $M$ be its exterior and $\tilde\mu$
its meridian.  If $K'$ is obtained from $K$ by the two-strand $n$-twist along $D$, then
\begin{equation}\label{eq:rogers}
\Sigma(K')=M(\tilde\mu'),\qquad \Delta(\tilde\mu,\tilde\mu')=|n|
\end{equation}
for a slope $\tilde\mu'$ on $\partial M$ \cite[Lemma~2.8]{Rog}, and again this identification respects orientations.

\begin{lemma}\label{lem:band}
Suppose that $K'$ is obtained from $K$ by a two-strand twist of order one and that $\det K'=\det K$.  Then
$\Delta(\tilde\mu,\lambda_M)=\Delta(\tilde\mu',\lambda_M)=1$.  Moreover, if $\tilde\mu$ and $\lambda_M$ are oriented
so that $\tilde\mu\cdot\lambda_M=1$ on $\partial\nu(\tg)$, then $\tilde\mu'=\pm(\tilde\mu+\delta\lambda_M)$ for some
$\delta\in\{1,-1\}$.
\end{lemma}

\begin{proof}
As in the proof of Lemma~\ref{lem:LM}, \eqref{eq:watson} gives
$\Delta(\tilde\mu,\lambda_M)=\Delta(\tilde\mu',\lambda_M)=m$, and writing $\tilde\mu=m\rho+r\lambda_M$ and
$\tilde\mu'=\pm(m\rho+s\lambda_M)$ with $\rho\cdot\lambda_M=1$, we obtain $1=\Delta(\tilde\mu,\tilde\mu')=m\,|r-s|$.
Hence $m=1$ and $|r-s|=1$, and the choice $\rho=\tilde\mu$ gives the second assertion.
\end{proof}

For twists of even order we use all the twists along the same twisting disk $D$, with the product neighborhood $B$
fixed above.  Then $B$ is a ball around $\gamma$ whose boundary meets $K$ in four points, so that $M$ is the double
cover of $S^3\setminus\operatorname{int}B$ branched along $K\setminus B$, and the two-strand $k$-twist along $D$ is
obtained by regluing $(B,B\cap K)$ along the $k$-th power of a half-twist along the copy $\partial D\times\{1/2\}$ of
$c$ in $\partial B$ \cite[Section~2.4]{Rog}.  The preimage of this copy
of $c$ in $\partial M$ consists of two parallel curves, whose slope we denote by $\ell$, and the half-twist lifts to a
Dehn twist along $\ell$.  Since $c$ meets a meridian arc of $K$ in $\partial B$ in one point, $\Delta(\tilde\mu,\ell)=1$
\cite[proof of Proposition~2.9]{Rog}, and we orient $\tilde\mu$ and $\ell$ so that $\tilde\mu\cdot\ell=1$.  Hence there
is a sign $\varepsilon\in\{1,-1\}$ such that, for every integer $k$, the knot $K_k$ obtained from $K$ by the two-strand
$k$-twist along $D$, where $K_0=K$, satisfies
\begin{equation}\label{eq:twistfamily}
\Sigma(K_k)=M(\tilde\mu+\varepsilon k\ell),
\end{equation}
and again these identifications respect orientations.  Each $K_k$ is a knot \cite[Remark~2.7]{Rog}.  Indeed, since
the two arcs of $K\cap B$ cross $D$ in opposite directions, the two arcs of $K\setminus\operatorname{int}B$ join the
two points of $K\cap(D\times\{0\})$ to each other and the two points of $K\cap(D\times\{1\})$ to each other, and
together with any two disjoint arcs in $B$, each joining a point of $D\times\{0\}$ to a point of $D\times\{1\}$, they
form a single circle.

\begin{lemma}\label{lem:gcc}
Let $n$ be a nonzero even integer, and suppose that $\det K_n=\det K$.  Then $\lambda_M=\ell$.  In particular,
$\Delta(\tilde\mu,\lambda_M)=1$, and if $\lambda_M$ is oriented so that $\tilde\mu\cdot\lambda_M=1$, then
$\Sigma(K_n)=M(\tilde\mu+\varepsilon n\lambda_M)$.
\end{lemma}

\begin{proof}
Write $\lambda_M=x\tilde\mu+y\ell$ with $\gcd(x,y)=1$.  For every integer $k$, the group
$H_1(M(\tilde\mu+\varepsilon k\ell))=H_1(\Sigma(K_k))$ is finite, since $K_k$ is a knot.  Hence
$\tilde\mu+\varepsilon k\ell\neq\lambda_M$, and~\eqref{eq:watson} gives $\det K_k=c_M\,|y-\varepsilon kx|$.  The
hypothesis $\det K_n=\det K_0$ therefore gives $|y-\varepsilon nx|=|y|$, so either $x=0$ or $2y=\varepsilon nx$.  In the
second case $x$ divides $y$, since $n$ is even, so $x=\pm1$ and $\lambda_M=\pm(\tilde\mu+\varepsilon(n/2)\ell)$, which is
impossible by the case $k=n/2\in\Z$, since $K_{n/2}$ is a knot.  Hence $x=0$, so $\lambda_M=\ell$ as slopes, and the remaining assertions follow from
$\tilde\mu\cdot\ell=1$.
\end{proof}

\subsection{Criteria for nugatory crossing changes}\label{subsec:nugatory}
A rational homology sphere $Y$ is an \emph{L-space} if $\HFhat(Y)$ is a free abelian group of rank $|H_1(Y)|$
\cite[Definition~1.1]{OSlens}.  Lens spaces are L-spaces \cite[Proposition~2.3]{OSlens}, and so are the branched
double covers of quasi-alternating links \cite[Proposition~3.3]{OS05}, which include all links with a connected
alternating diagram \cite[Lemma~3.2]{OS05}.  The Dehn surgery characterization of the unknot is due to Kronheimer,
Mrowka, Ozsv\'ath and Szab\'o for $Y=S^3$ \cite[Theorem~1.1]{KMOS}, and to Gainullin for null-homologous knots in
L-spaces \cite[Theorem~8.2]{Gai}.  We state it
in the form of \cite[Theorem~3.1]{LM}.  Here an \emph{unknot} $U\subset Y$ is a knot that bounds an embedded disk,
and surgery coefficients of null-homologous knots refer to their Seifert framings.

\begin{theorem}[Kronheimer--Mrowka--Ozsv\'ath--Szab\'o, Gainullin]\label{thm:gainullin}
Let $\kappa$ be a null-homologous knot in an L-space $Y$, and let $U\subset Y$ be an unknot.  If
$Y_{p/q}(\kappa)\cong Y_{p/q}(U)$ for some nonzero $p/q\in\Q$, then $\kappa$ is isotopic to $U$.
\end{theorem}

Gainullin works with coefficients in an arbitrary field \cite[Section~2]{Gai}, so, as in \cite[Corollary~1.3]{LM},
Theorem~\ref{thm:gainullin} holds more generally when $Y$ is an L-space with coefficients in a field, for instance $\Z/2$, that is, when
$\HFhat(Y;\Z/2)$ has rank $|H_1(Y)|$.  Since L-spaces enter the results below only through
Theorem~\ref{thm:gainullin}, the same applies to them.

The second ingredient is the equivariant Dehn lemma, in the form of the following result of Lidman and Moore
\cite[Proposition~3.3]{LM}.  It is stated there under the standing assumption that the crossing change preserves the
isotopy type of the knot, and its proof uses this assumption only through the conclusion of Lemma~\ref{lem:LM}.

\begin{proposition}[Lidman--Moore]\label{prop:LMnugatory}
In the setting of Section~\ref{subsec:montesinos}, suppose that $\det K_+=\det K_-$ and that $\tg$ bounds an embedded
disk in $\Sigma(K_+)$.  Then the crossing change at $c$ is nugatory.
\end{proposition}

\begin{corollary}[Lidman--Moore]\label{cor:LMnugatory}
In the setting of Section~\ref{subsec:montesinos}, suppose that $\Sigma(K_+)$ is an L-space, that
$\Sigma(K_-)\cong\Sigma(K_+)$, and that $\tg$ is null-homologous.  Then the crossing change at $c$ is nugatory.
\end{corollary}

\begin{proof}
Since $\det K_\pm=|H_1(\Sigma(K_\pm))|$, Lemma~\ref{lem:LM} shows that $\Delta(\alpha,\lambda_M)=1$ and
$\beta=\pm(\alpha+2\varepsilon\lambda_M)$.  Since $\tg$ is null-homologous, $\lambda_M$ is its Seifert framing, so
$\Sigma(K_-)=M(\beta)=Y_{\varepsilon/2}(\tg)$, where $Y=\Sigma(K_+)$.  For an unknot $U\subset Y$ we have
$Y_{\varepsilon/2}(U)\cong Y$, so $Y_{\varepsilon/2}(\tg)\cong Y_{\varepsilon/2}(U)$, and Theorem~\ref{thm:gainullin}
shows that $\tg$ is isotopic to $U$.  By Proposition~\ref{prop:LMnugatory}, the crossing change is nugatory.
\end{proof}

This is the observation of \cite[Remark~3.4]{LM}: if $\Sigma(K)$ is an L-space and the lift of a crossing arc is
null-homologous, then the crossing change is nugatory or changes the isotopy type of $K$.

Rogers extended Proposition~\ref{prop:LMnugatory} to two-strand twists of order at least two.  We need it for
generalized crossing changes.

\begin{proposition}[Rogers]\label{prop:rogersgcc}
Let $n$ be a nonzero even integer, and let $K'$ be obtained from a knot $K\subset S^3$ by the two-strand $n$-twist
along a twisting disk $D$, with twisting arc $\gamma\subset D$.  If $\det K'=\det K$ and $\tg$ bounds an embedded disk
in $\Sigma(K)$, then the twist is nugatory.
\end{proposition}

\begin{proof}
This is \cite[Proposition~2.9]{Rog}, where $K'$ is assumed to be isotopic to $K$.  Its proof uses this assumption only
to show that $\Delta(\tilde\mu',\lambda)=1$, where $\Sigma(K')=M(\tilde\mu')$ and $\lambda$ is the slope on $\partial M$
that bounds a disk in $M$.  This slope is the rational longitude of $M$, so $\Delta(\tilde\mu',\lambda)=1$ follows from
$\det K'=\det K$ by Lemma~\ref{lem:gcc}.
\end{proof}

\begin{corollary}\label{cor:rogersgcc}
In the setting of Proposition~\ref{prop:rogersgcc}, suppose that $\Sigma(K)$ is an L-space, that
$\Sigma(K')\cong\Sigma(K)$, and that $\tg$ is null-homologous.  Then the twist is nugatory.
\end{corollary}

\begin{proof}
Since $\det K'=\det K$, Lemma~\ref{lem:gcc} shows that $\Sigma(K')=M(\tilde\mu+\varepsilon n\lambda_M)$.  Since $\tg$
is null-homologous, $\lambda_M$ is its Seifert framing, so $\Sigma(K')=Y_{\varepsilon/n}(\tg)$, where $Y=\Sigma(K)$, and
$Y_{\varepsilon/n}(U)\cong Y$ for an unknot $U\subset Y$.  By Theorem~\ref{thm:gainullin}, $\tg$ is isotopic to $U$,
and Proposition~\ref{prop:rogersgcc} applies.
\end{proof}

The analogues of Propositions~\ref{prop:LMnugatory} and~\ref{prop:rogersgcc} fail for two-strand twists of order one.
Rogers exhibited a
two-strand twist of order one on the unknot that produces an unknot and is not nugatory \cite[Example~4.1]{Rog}.  He
introduced the following weaker notion \cite[Definition~4.4]{Rog}.  A two-strand twist of order one along $D$, with
product neighborhood $B$, that produces $K'$ is \emph{weakly nugatory} if it is nugatory, or if there is a twisting
circle $c'\subset\partial B$ for $K$ that bounds a disk in $S^3\setminus K$, meets the copy of $c$ in $\partial B$ in
two points and a meridian arc of $K$ in $\partial B$ in one point, and such that regluing $(B,B\cap K)$ along a
half-twist along $c'$ of the opposite sign also produces $K'$.  In either case, some nugatory two-strand twist of order
one on $K$ produces a knot isotopic to $K'$: the given twist in the first case, and the twist along $c'$ in the second.
By Remark~\ref{rem:bandmutants}, this does not imply that $K'$ is isotopic to $K$.

\begin{proposition}[Rogers]\label{prop:rogers}
Let $K'$ be obtained from a knot $K\subset S^3$ by a two-strand twist of order one with twisting arc $\gamma$.  If
$\det K'=\det K$ and $\tg$ bounds an embedded disk in $\Sigma(K)$, then the twist is weakly nugatory.
\end{proposition}

\begin{proof}
This is \cite[Proposition~4.6]{Rog} under a weaker hypothesis.  As in the proof of
Proposition~\ref{prop:rogersgcc}, the assumption made there that $K'$ is isotopic to $K$ is used only to show that $\Delta(\tilde\mu',\lambda)=1$, where $\lambda$ is the
slope on $\partial M$ that bounds a disk in $M$ \cite[proof of Proposition~2.9]{Rog}.  This slope is the rational longitude of $M$, so
$\Delta(\tilde\mu',\lambda)=1$ follows from $\det K'=\det K$ by Lemma~\ref{lem:band}.
\end{proof}

\begin{corollary}\label{cor:rogers}
Let $K'$ be obtained from a knot $K\subset S^3$ by a two-strand twist of order one with twisting arc $\gamma$.  If
$\Sigma(K)$ is an L-space, $\Sigma(K')\cong\Sigma(K)$, and $\tg$ is null-homologous, then the twist is weakly
nugatory.
\end{corollary}

\begin{proof}
By Lemma~\ref{lem:band}, $\Delta(\tilde\mu,\lambda_M)=1$ and $\tilde\mu'=\pm(\tilde\mu+\delta\lambda_M)$.  Since
$\tg$ is null-homologous, $\lambda_M$ is its Seifert framing, so $\Sigma(K')=Y_\delta(\tg)$ with $Y=\Sigma(K)$, and
$Y_\delta(U)\cong Y$ for an unknot $U\subset Y$.  By Theorem~\ref{thm:gainullin}, $\tg$ is isotopic to $U$, and
Proposition~\ref{prop:rogers} applies.
\end{proof}

\section{Monotonicity of correction terms}\label{sec:monotone}

Throughout this section, $Y$ is a rational homology sphere and $\kappa\subset Y$ is an isotropic knot, with exterior
$M$ and meridian $\mu$, and $\lambda=\lambda_M$ is the rational longitude of $M$.  By Lemma~\ref{lem:isotropic}(a),
$\lambda$ is a framing of $\kappa$, and we orient $\mu$ and $\lambda$ so that $\mu\cdot\lambda=1$ on
$\partial\nu(\kappa)$.  We denote by $o$ the order of $[\kappa]$ in $H_1(Y)$.  For $j\in\Z$ we put
\begin{equation}\label{eq:family}
\mu_j=\mu-j\lambda\qquad\text{and}\qquad Z_j=M(\mu_j),
\end{equation}
and we denote by $\kappa_j\subset Z_j$ the dual knot.  Thus $Z_0=Y$ and $\kappa_0=\kappa$, and for $j\neq0$ the
manifold $Z_j$ is the result of surgery on $\kappa$ with coefficient $-1/j$ with respect to the framing $\lambda$.
The main result of this section is the following restatement of Theorem~\ref{thm:C}.

\begin{theorem}\label{thm:monotone}
For every $j\in\Z$,
\begin{equation}\label{eq:monotone}
S(Z_{j+1})-S(Z_j)\ \ge\ \frac{|H_1(Y)|\,(o^2-1)}{6\,o^2}.
\end{equation}
In particular, if $\kappa$ is not null-homologous, then $Z_j\not\cong Z_{j'}$ for all $j\neq j'$, and
$M(\eta)\not\cong Y$ for every slope $\eta\neq\mu$.
\end{theorem}

For a null-homologous knot the bound~\eqref{eq:monotone} is zero.  If $Y$ is an integral homology sphere, then
$S(Y)=d(Y)$, and in this case Theorem~\ref{thm:monotone} gives the inequalities
\[
d(Y_{1/(n+1)}(\kappa))\le d(Y_{1/n}(\kappa))\le d(Y)\le d(Y_{-1/n}(\kappa))\le d(Y_{-1/(n+1)}(\kappa))\qquad(n\ge1)
\]
of \cite[Corollary~9.14]{OS03}.  The bound~\eqref{eq:monotone} is sharp; see Section~\ref{subsec:examples}.

\subsection{The two-handle cobordism}\label{subsec:cobordism}
We first show that each step from $Z_j$ to $Z_{j+1}$ is a surgery of the same kind as the step from $Z_0$ to $Z_1$.

\begin{lemma}\label{lem:family}
Let $j\in\Z$.  The knot $\kappa_j\subset Z_j$ has exterior $M$, meridian $\mu_j$ and rational longitude $\lambda$, and
$\mu_j\cdot\lambda=1$.  It is isotropic, its homology class has order $o$, and $|H_1(Z_j)|=|H_1(Y)|$.  Moreover,
$Z_{j+1}$ is obtained from $Z_j$ by surgery on $\kappa_j$ along the framing $\mu_j-\lambda$, which has coefficient
$-1$ with respect to $\lambda$.
\end{lemma}

\begin{proof}
The exterior and the meridian of $\kappa_j$ are $M$ and $\mu_j$ by construction, and the rational longitude depends
only on $M$.  The tori $\partial\nu(\kappa_j)$ and $\partial\nu(\kappa)$ are both identified with $\partial M$, with
the same orientation, so $\mu_j\cdot\lambda=\mu\cdot\lambda=1$.  In particular $\Delta(\mu_j,\lambda)=1$, and the
second assertion follows from Lemma~\ref{lem:isotropic}(c).  Finally,
$\Delta(\mu_j,\mu_j-\lambda)=|\mu_j\cdot\lambda|=1$,
so $\mu_j-\lambda$ is a framing of $\kappa_j$; it has coefficient $1/(-1)=-1$ with respect to $\lambda$, and
$M(\mu_j-\lambda)=M(\mu_{j+1})=Z_{j+1}$.
\end{proof}

By Lemma~\ref{lem:family}, it suffices to analyze the following situation, which we fix for the rest of this
subsection and the next one.  Let $Y'=M(\mu-\lambda)$, let $\kappa'\subset Y'$ be its dual knot, and let $W$ be the
cobordism from $Y$ to $Y'$ obtained from $Y\times[0,1]$ by attaching a $2$-handle along $\kappa\times\{1\}$ with
framing $\mu-\lambda$.  Let $D\in H_2(W,Y)$ be the class of the core of the $2$-handle, extended by the annulus
$\kappa\times[0,1]$, and let $D^*\in H_2(W,Y')$ be the class of the cocore, whose boundary is isotopic to $\kappa'$ in
$Y'$.  The core and the cocore meet transversely in one point, so $D\cdot D^*=\pm1$ for the intersection pairing
$H_2(W,Y)\times H_2(W,Y')\to\Z$.

\begin{lemma}\label{lem:W}
\begin{enumerate}[label=(\alph*)]
\item $H_1(W)\cong H_1(Y)/\langle[\kappa]\rangle\cong H_1(Y')/\langle[\kappa']\rangle$, moreover  $b_1(W)=0$.
\item $H_2(W)\cong\Z$, generated by a class $A$ whose images in $H_2(W,Y)$ and in $H_2(W,Y')$ are $\pm o\,D$ and
$\pm o\,D^*$, respectively.
\item $A\cdot A=-o^2$.  In particular, $b_2^+(W)=0$.
\item The restriction maps $H^2(W)\to H^2(Y)$ and $H^2(W)\to H^2(Y')$ are surjective, and their kernels are infinite
cyclic, generated by classes $\theta$ and $\theta'$ with $\langle\theta,A\rangle=\pm o$ and $\langle\theta',A\rangle=\pm
o$.
\end{enumerate}
\end{lemma}

\begin{proof}
Relative to $Y$, the cobordism $W$ consists of a single $2$-handle, so $H_*(W,Y)$ is infinite cyclic in degree two,
generated by $D$, and vanishes in all other degrees.  Turned upside down, $W$ is obtained from $Y'\times[0,1]$ by
attaching a single $2$-handle along $\kappa'$, whose core is the cocore of the original one; hence $H_*(W,Y')$ is
infinite cyclic in degree two, generated by $D^*$, and vanishes in all other degrees.

(a), (b) Since $Y$ is a rational homology sphere, $H_2(Y)=0$, and the exact sequence of the pair $(W,Y)$ reads
\[
0\longrightarrow H_2(W)\longrightarrow H_2(W,Y)=\Z D\xrightarrow{\ \partial\ }H_1(Y)\longrightarrow
H_1(W)\longrightarrow0,
\]
with $\partial D=[\kappa]$.  It gives the first isomorphism of~(a), and it identifies $H_2(W)$ with the subgroup of
$\Z D$ generated by $o\,D$; let $A$ be the class corresponding to $o\,D$.  Since $Y'$ is a rational homology sphere by
Lemma~\ref{lem:isotropic}(c), the exact sequence of the pair $(W,Y')$, in which $\partial D^*=[\kappa']$, likewise
gives the second isomorphism of~(a), and it identifies $H_2(W)$ with the subgroup of $\Z D^*$ generated by $o'D^*$,
where $o'$ is the order of $[\kappa']$ in $H_1(Y')$.  By Lemma~\ref{lem:isotropic}(c), applied to $\eta=\mu-\lambda$,
we have $o'=o$.  This proves~(b).

(c) Present $Y=\partial X_N$ as in Section~\ref{subsec:isotropic}, with $\kappa\subset S^3\setminus L$ and linking
vector $v$.  By Lemma~\ref{lem:presentation}, $x_0=v^TN^{-1}v$ is an integer and $\lambda$ has coefficient $x_0$ with
respect to the Seifert framing of $\kappa$ in $S^3$, so the framing $\mu-\lambda$ has coefficient $x_0-1$.  Hence
$X_B=X_N\cup_YW$ is the $2$-handlebody of the framed link $L\cup\kappa$, with linking matrix
\[
B=\begin{pmatrix}N&v\\ v^T&x_0-1\end{pmatrix}.
\]
Let $e_1,\dots,e_n,g$ be the basis of $H_2(X_B)$ given by the $2$-handles, with $g$ attached along $\kappa$, put
$\hat v=N^{-1}v\in\Q^n$, and let
\[
g'=g-\sum_{i=1}^n\hat v_ie_i\in H_2(X_B;\Q).
\]
Then $g'\cdot e_i=0$ for all $i$, and $g'\cdot g'=x_0-1-2\hat v^Tv+\hat v^TN\hat v=-1$.  Since $\det B=-\det N\neq0$,
the orthogonal
complement of $H_2(X_N;\Q)$ in $H_2(X_B;\Q)$ is the line $\Q g'$.  Classes in $H_2(W)$ and in $H_2(X_N)$ are
represented by disjoint cycles, so the image of $H_2(W)$ in $H_2(X_B)$ is contained in $\Q g'\cap H_2(X_B)$, which is
$\Z\,og'$ because the least common multiple of the denominators of the entries of $\hat v$ is $o$.  On the other hand,
by
excision and Lefschetz duality,
\begin{gather*}
H_3(X_B,W)\cong H_3(X_N,Y)\cong H^1(X_N)=0,\\
H_2(X_B,W)\cong H_2(X_N,Y)\cong H^2(X_N)\cong\Z^n,
\end{gather*}
so the map $H_2(W)\to H_2(X_B)$ is injective with torsion-free cokernel.  It follows that $H_2(W)$ maps onto $\Z\,og'$,
so $A=\pm og'$ and $A\cdot A=o^2\,(g'\cdot g')=-o^2$.

(d) Consider the exact sequence
\[
H^1(Y)\longrightarrow H^2(W,Y)\longrightarrow H^2(W)\longrightarrow H^2(Y)\longrightarrow H^3(W,Y).
\]
By Lefschetz duality, $H^3(W,Y)\cong H_1(W,Y')=0$ and $H^2(W,Y)\cong H_2(W,Y')=\Z D^*$, and $H^1(Y)=0$.  Hence the
restriction to $Y$ is surjective, and its kernel is infinite cyclic, generated by the image $\theta$ of the class
dual to $D^*$.  By~(b),
\[
\langle\theta,A\rangle=\pm\,D^*\cdot A=\pm\,D^*\cdot(o\,D)=\pm o .
\]
The same argument at $Y'$, with the roles of $D$ and $D^*$ exchanged, gives a generator $\theta'$ of the kernel of the
restriction to $Y'$ with $\langle\theta',A\rangle=\pm\,D\cdot(o\,D^*)=\pm o$.
\end{proof}

The cobordism $W$ was also considered by Caudell \cite{Cau}, who showed that, for an integral surgery between rational
homology spheres, the generator of the second homology of the trace has self-intersection of absolute value
$o\,o'$, where $o$ and $o'$ are the orders of the classes of the knot and of its dual \cite[Proposition~6]{Cau}.
Whenever $o\ge2$, the second assertion of \cite[Proposition~6]{Cau}, that $o$ and $o'$ are coprime, is incompatible
with Lemma~\ref{lem:W}(b), which gives $o'=o$ for the cobordisms considered here; see Examples~\ref{ex:pretzel}
and~\ref{ex:BHW}.

\subsection{Spin\texorpdfstring{$^c$}{c} structures}\label{subsec:spinc}
For $\mathfrak s\in\Spinc(W)$ we put
\[
c_0(\mathfrak s)=\langle c_1(\mathfrak s),A\rangle\in\Z .
\]

\begin{lemma}\label{lem:fibres}
\begin{enumerate}[label=(\alph*)]
\item For every $\mathfrak s\in\Spinc(W)$, we have $c_0(\mathfrak s)\equiv o\pmod 2$ and $c_1(\mathfrak
s)^2=-c_0(\mathfrak s)^2/o^2$.
\item For every $\mathfrak t\in\Spinc(Y)$, the set of $\mathfrak s\in\Spinc(W)$ with $\mathfrak s|_Y=\mathfrak t$ is
nonempty, and $c_0$ maps it bijectively onto a coset of $2o\Z$ in $\Z$.  The same holds for $\mathfrak t'\in\Spinc(Y')$
and the restriction to $Y'$.
\end{enumerate}
\end{lemma}

\begin{proof}
(a) The class $c_1(\mathfrak s)$ is characteristic, so $c_0(\mathfrak s)\equiv A\cdot A=-o^2\equiv o\pmod2$.  Since
$H_2(W;\Q)=\Q A$, the image of $c_1(\mathfrak s)$ in $H^2(W;\Q)$ is the restriction of the class in
$H^2(W,\partial W;\Q)$ that is Lefschetz dual to $(c_0(\mathfrak s)/A\cdot A)\,A$, so $c_1(\mathfrak
s)^2=c_0(\mathfrak s)^2/(A\cdot A)=-c_0(\mathfrak s)^2/o^2$.

(b) The set $\Spinc(W)$ is an affine space over $H^2(W)$, the restriction map to $\Spinc(Y)$ is equivariant with
respect to the restriction $H^2(W)\to H^2(Y)$, and $c_1(\mathfrak s+x)=c_1(\mathfrak s)+2x$.  By
Lemma~\ref{lem:W}(d), the restriction map is therefore surjective, and its fibres are the orbits of the free action
of $\Z\theta$.  Since $c_0(\mathfrak s+k\theta)=c_0(\mathfrak s)\pm2ko$, the assertion follows.  The same argument
applies at $Y'$, with $\theta'$ in place of $\theta$: the two kernels may be different subgroups of $H^2(W)$, but
$\langle\theta',A\rangle=\pm o$ by Lemma~\ref{lem:W}(d), so $c_0$ again changes by $\pm2o$ along each fibre.
\end{proof}

Let
\[
\mathcal S_W=\{\mathfrak s\in\Spinc(W):-o<c_0(\mathfrak s)\le o\}.
\]

\begin{lemma}\label{lem:window}
\begin{enumerate}[label=(\alph*)]
\item The restriction maps $\mathcal S_W\to\Spinc(Y)$ and $\mathcal S_W\to\Spinc(Y')$ are bijective.
\item For every integer $a$ with $-o<a\le o$ and $a\equiv o\pmod 2$, exactly $|H_1(Y)|/o$ elements $\mathfrak
s\in\mathcal S_W$ satisfy $c_0(\mathfrak s)=a$.
\end{enumerate}
\end{lemma}

\begin{proof}
(a) The interval $(-o,o]$ contains exactly one integer from each residue class modulo $2o$, so this follows from
Lemma~\ref{lem:fibres}(b).

(b) For $\mathfrak t\in\Spinc(Y)$, let $r(\mathfrak t)\in\Z/2o$ be the residue of $c_0$ on the fiber over $\mathfrak
t$.  For $x\in H^2(Y)$, choose $\tilde x\in H^2(W)$ with $\tilde x|_Y=x$, and let $\varphi(x)\in\Z/o$ be the residue
of $\langle\tilde x,A\rangle$.  By Lemma~\ref{lem:W}(d), two choices of $\tilde x$ differ by a multiple of $\theta$,
so $\varphi\colon H^2(Y)\to\Z/o$ is a well-defined homomorphism.  It is surjective.  Indeed, by the universal coefficient
theorem, evaluation on $A$ maps $H^2(W)$ onto $\operatorname{Hom}(H_2(W),\Z)\cong\Z$, so there is $\tilde y\in
H^2(W)$ with $\langle\tilde y,A\rangle=1$, and then $\varphi(\tilde y|_Y)=1$.  If $\mathfrak s|_Y=\mathfrak t$, then
$(\mathfrak s+\tilde x)|_Y=\mathfrak t+x$ and $c_0(\mathfrak s+\tilde x)=c_0(\mathfrak s)+2\langle\tilde x,A\rangle$,
so
\[
r(\mathfrak t+x)=r(\mathfrak t)+2\varphi(x)\in\Z/2o ,
\]
where $2\varphi(x)$ denotes the image of $\varphi(x)$ under the injection $\Z/o\to\Z/2o$ given by multiplication by
two.  Fix $\mathfrak t_0\in\Spinc(Y)$.  As $x$ runs over $H^2(Y)$, which has $|H_1(Y)|$ elements, $2\varphi(x)$ takes
each value in $2\Z/2o\Z$ exactly $|H_1(Y)|/o$ times.  Since $r(\mathfrak t_0)\equiv o\pmod2$ by
Lemma~\ref{lem:fibres}(a), it follows that $r$ takes each of the $o$ residues modulo $2o$ that are congruent to $o$
modulo $2$ exactly $|H_1(Y)|/o$ times.  By~(a), each $\mathfrak t$ is the restriction of exactly one $\mathfrak
s\in\mathcal S_W$, and $c_0(\mathfrak s)$ is the representative of $r(\mathfrak t)$ in $(-o,o]$; this proves~(b).
\end{proof}

\subsection{Proof of Theorem~\ref{thm:monotone}}\label{subsec:proofmono}
Fix $j\in\Z$.  By Lemma~\ref{lem:family}, the knot $\kappa_j\subset Z_j$ is isotropic of order $o$, with meridian
$\mu_j$ and rational longitude $\lambda$, and $Z_{j+1}=M(\mu_j-\lambda)$.  We apply the construction of
Section~\ref{subsec:cobordism} to $Z_j$ and $\kappa_j$ in place of $Y$ and $\kappa$, so that $Y'=Z_{j+1}$.  By
Lemma~\ref{lem:W}(c), $b_2^+(W)=0$, so Proposition~\ref{prop:dineq} applies.  By Lemma~\ref{lem:window}(a), every
spin$^c$ structure on $Z_j$ and on $Z_{j+1}$ is the restriction of exactly one element of $\mathcal S_W$.  Summing
\eqref{eq:dineq} over $\mathcal S_W$ and using Lemma~\ref{lem:fibres}(a), we obtain
\[
S(Z_{j+1})-S(Z_j)=\sum_{\mathfrak s\in\mathcal S_W}\bigl(d(Z_{j+1},\mathfrak s|_{Z_{j+1}})-d(Z_j,\mathfrak
s|_{Z_j})\bigr)\ \ge\ \sum_{\mathfrak s\in\mathcal S_W}\frac{o^2-c_0(\mathfrak s)^2}{4o^2}.
\]
By Lemma~\ref{lem:window}(b) and since $|H_1(Z_j)|=|H_1(Y)|$, the integers $c_0(\mathfrak s)$ with $\mathfrak
s\in\mathcal S_W$ are the numbers $o-2k$ with $0\le k\le o-1$, each taken $|H_1(Y)|/o$ times.  Since
\[
\sum_{k=0}^{o-1}\bigl(o^2-(o-2k)^2\bigr)=\sum_{k=0}^{o-1}(4ok-4k^2)=2o^2(o-1)-\frac{2(o-1)o(2o-1)}{3}
=\frac{2o(o^2-1)}{3},
\]
the right-hand side equals $|H_1(Y)|(o^2-1)/(6o^2)$, which proves~\eqref{eq:monotone}.

Suppose now that $\kappa$ is not null-homologous, so that $o\ge2$ and $S(Z_j)$ is strictly increasing in $j$.  Since
$S$ depends only on the orientation-preserving homeomorphism type, $Z_j\not\cong Z_{j'}$ for $j\neq j'$.  If $\eta$ is a
slope with $M(\eta)\cong Y$, then $\eta\neq\lambda$, since $H_1(M(\lambda))$ is infinite, and
$|H_1(M(\eta))|=|H_1(Y)|$, so $\Delta(\eta,\lambda)=\Delta(\mu,\lambda)=1$ by~\eqref{eq:watson}.  Writing
$\eta=a\mu+b\lambda$, we get $|a|=|\eta\cdot\lambda|=1$, so $\eta=\pm\mu_j$ with $j=-ab$.  Thus $M(\eta)=Z_j\cong Z_0$,
and therefore $j=0$ and $\eta=\mu$.\qed

\subsection{Examples}\label{subsec:examples}
The following examples show that the inequality~\eqref{eq:monotone} is attained, and that the last assertion of
Theorem~\ref{thm:monotone} does not extend to orientation-reversing homeomorphisms; Corollary~\ref{cor:reverse} below
shows what the inequality still says about them.  In both of them, the correction terms of lens spaces are computed
with the recursive formula of \cite[Proposition~4.8]{OS03}; one finds
\begin{equation}\label{eq:lensS}
|S(L(9,2))|=\frac43\qquad\text{and}\qquad |S(L(49,18))|=4 .
\end{equation}

\begin{example}\label{ex:pretzel}
For $n\in\Z$, let $K_n=P(3,-3,n)$ be the pretzel knot whose three twist regions contain $3$, $-3$ and $n$ half-twists.
Then $\det K_n=|3\cdot(-3)+(-3)\,n+n\cdot3|=9$ for every $n$.  The knot $K_0$ is the connected sum of the two
trefoils, $K_1$ is the stevedore knot $6_1$, which is the two-bridge knot $b(9,2)$, up to mirror image, and $K_{-1}$
is the mirror image of $K_1$, since a pretzel knot with three twist regions does not depend on their order.  The
knots $K_n$ agree outside a ball $\Omega$ that meets each of them in the third twist region, an integer tangle whose
two strands are oppositely oriented.  Let $M$ be the double cover of $S^3\setminus\Omega$ branched along
$K_n\setminus\Omega$, which does not depend on $n$, and let $\eta_n$ be the slope on $\partial M$ with
$\Sigma(K_n)=M(\eta_n)$, namely the slope of the two parallel curves that form the preimage of the boundary of a disk in $\Omega$
that separates the two strands of $K_n$.
For each $n$, the ball $\Omega$ contains a twisting disk for $K_n$ whose two-strand $k$-twists turn $K_n$ into
$K_{n+k}$, and the complement in $\Omega$ of a small ball around this disk meets $K_n$ in a product $P\times[0,1]$,
where $P$ consists of four points, so the double cover of this complement branched along $P\times[0,1]$ is a collar of
$\partial M$.  Hence $M$ is identified with the exterior of the lift of the corresponding twisting arc,
and~\eqref{eq:rogers} gives $\Delta(\eta_n,\eta_{n+1})=1$ and
$\Delta(\eta_{n-1},\eta_{n+1})=2$.  By Lemma~\ref{lem:band}, $\Delta(\eta_n,\lambda_M)=1$ for all $n$.

Let $Y=\Sigma(K_{-1})$, let $\kappa\subset Y$ be the dual knot of $M(\eta_{-1})$, which is the lift of the arc of a
crossing change in $\Omega$, and define $Z_j$ by~\eqref{eq:family} with $\mu=\eta_{-1}$.  The slopes at distance one
from $\lambda_M$ are the slopes $\mu_j$, and $\Delta(\mu_j,\mu_{j'})=|j-j'|$; hence there is a sign $\tau\in\{1,-1\}$
such that $\eta_n=\pm\mu_{\tau(n+1)}$ and $\Sigma(K_n)=Z_{\tau(n+1)}$ for all $n$.  The class $[\kappa]$ has order
$3$.  Indeed, $H_1(Y)\cong\Z/9$, since $Y$ is the lens space $\pm L(9,2)$, so the isotropic classes have order one or
three.  If $\kappa$ were null-homologous, then by Lemma~\ref{lem:W}(a), applied to the cobordism between $Y$ and
$Z_{\pm1}=\Sigma(K_0)$, the group $H_1(\Sigma(K_0))\cong(\Z/3)^2$ would have a quotient isomorphic to $H_1(Y)\cong\Z/9$.
By Lemma~\ref{lem:family}, the dual knot in $\Sigma(K_0)=Z_\tau$ of the integral surgery from $Y$ to $\Sigma(K_0)$ also
has order three.

Thus $o=3$ and $|H_1(Y)|=9$, and~\eqref{eq:monotone} reads $S(Z_{j+1})-S(Z_j)\ge4/3$.  On the other hand,
$\Sigma(K_1)\cong-\Sigma(K_{-1})$, and $S(\Sigma(K_0))=S(L(3,1)\#{-L(3,1)})=0$ by the additivity of correction terms
under connected sums \cite[Theorem~4.3]{OS03}.  By~\eqref{eq:lensS}, the two steps from $\Sigma(K_{-1})$ through
$\Sigma(K_0)$ to $\Sigma(K_1)$ change $S$ by exactly $4/3$ each, so equality holds in~\eqref{eq:monotone}.  The
crossing change in $\Omega$ that turns $K_{-1}$ into $K_1$ changes $6_1$ into its mirror image; the lift $\kappa$ of
its crossing arc is not null-homologous, and $M(\eta_1)=Z_{2\tau}\cong-Z_0$.
\end{example}

\begin{example}\label{ex:BHW}
Bleiler, Hodgson and Weeks described a knot $J\subset S^2\times S^1$ with hyperbolic exterior $X$, representing
seven times a generator of $H_1(S^2\times S^1)\cong\Z$, such that the fillings $X_1$ and $X_2$ of $X$ along the slopes
$18/1$ and $19/1$ in the coordinates of \cite[Section~3]{BHW} are oppositely oriented copies of $L(49,18)$.  These two
slopes have distance one from each other and from the meridian $1/0$ of $J$.  Since filling along
the meridian gives $S^2\times S^1$, the meridian of $J$ is the rational longitude $\lambda$ of $X$.  Regarding
$S^2\times S^1$ as $0$-surgery on an unknot $U\subset S^3$ with $\operatorname{lk}(J,U)=7$, we find that
$H_1(X)\cong\Z\oplus\Z/7$, where $i_*(\lambda)$ generates the torsion subgroup.  Hence, by
Lemma~\ref{lem:isotropic}(a) and~(b), the dual knot of $X_1$ is isotropic of order $o=7$, and $X_2=Z_{\pm1}$ in the
notation~\eqref{eq:family}, with $Y=X_1$.  Theorem~\ref{thm:monotone} gives $|S(X_2)-S(X_1)|\ge49\cdot48/(6\cdot49)=8$,
and by~\eqref{eq:lensS} equality holds.
\end{example}

In particular, the conclusion $M(\eta)\not\cong Y$ of Theorem~\ref{thm:monotone} does not extend to
orientation-reversing homeomorphisms.  In Examples~\ref{ex:pretzel} and~\ref{ex:BHW}, a nontrivial filling of the
exterior of a non-null-homologous isotropic knot in $Y$ yields $-Y$.  The inequality~\eqref{eq:monotone} still
restricts such fillings.

\begin{corollary}\label{cor:reverse}
In the setting of Theorem~\ref{thm:monotone}, suppose that $\kappa$ is not null-homologous, and let
$c=|H_1(Y)|(o^2-1)/(6o^2)$.  Then $Z_j\cong-Y$ for at most one $j\neq0$, and in that case $|j|\,c\le2|S(Y)|$.  In
particular, if $S(Y)=0$, then no nontrivial surgery on $\kappa$ yields $Y$ or $-Y$.
\end{corollary}

\begin{proof}
Since $S(-Y)=-S(Y)$, a homeomorphism $Z_j\cong-Y$ gives $S(Z_j)=-S(Y)$.  Since $o\ge2$, we have $c>0$, so
by~\eqref{eq:monotone} the sum $S(Z_j)$ is strictly increasing in $j$ and $|S(Z_j)-S(Z_0)|\ge|j|\,c$, where $Z_0=Y$.
This gives the first two assertions, and if $S(Y)=0$, then $j=0$ is the only solution of $S(Z_j)=-S(Y)$.  As in the
proof of Theorem~\ref{thm:monotone}, every slope $\eta$ with $|H_1(M(\eta))|=|H_1(Y)|$ is $\pm\mu_j$ for some $j$, so
every surgery on $\kappa$ that yields $Y$ or $-Y$ yields one of the manifolds $Z_j$.
\end{proof}

Both examples attain equality in Corollary~\ref{cor:reverse}: in Example~\ref{ex:pretzel}, $|j|=2$, $c=4/3$ and
$|S(Y)|=4/3$, and in Example~\ref{ex:BHW}, $|j|=1$, $c=8$ and $|S(Y)|=4$.

\section{Applications}\label{sec:applications}

\subsection{Cosmetic crossing changes}\label{subsec:ccc}
We use the notation of Section~\ref{subsec:montesinos}: $c$ is a crossing of a diagram of a knot, $K_+$ and $K_-$ are
the knots in which $c$ is positive and negative, $\gamma$ is a crossing arc, and $\tg\subset\Sigma(K_+)$ is its lift.
The arc $\gamma$ also joins the two points of $K_-\cap D$, and its lift to $\Sigma(K_-)=M(\beta)$ is the dual knot of
this filling.

\begin{theorem}\label{thm:ccc}
Suppose that $\Sigma(K_-)\cong\Sigma(K_+)$; for instance, $K_-$ and $K_+$ are isotopic, as oriented or as unoriented
knots.
\begin{enumerate}[label=(\alph*)]
\item The lifts of the crossing arc to $\Sigma(K_+)$ and to $\Sigma(K_-)$ are null-homologous.
\item If $\Sigma(K_+)$ is an L-space, then the crossing change at $c$ is nugatory; in particular, $K_+$ and $K_-$ are
isotopic.
\end{enumerate}
\end{theorem}

\begin{proof}
Let $Y=\Sigma(K_+)$, and let $M$, $\alpha$ and $\beta$ be as in Section~\ref{subsec:montesinos}.  Since
$|H_1(\Sigma(K_\pm))|=\det K_\pm$, the hypothesis gives $\det K_+=\det K_-$, so by Lemma~\ref{lem:LM},
$\Delta(\alpha,\lambda_M)=1$ and $\beta=\pm(\alpha+2\varepsilon\lambda_M)$ for some $\varepsilon\in\{1,-1\}$.  By
Lemma~\ref{lem:isotropic}(a), the knot $\tg$ is isotropic.  In the notation~\eqref{eq:family}, with $\kappa=\tg$,
$\mu=\alpha$ and $\lambda=\lambda_M$, we have $M(\beta)=M(\mu_{-2\varepsilon})=Z_{-2\varepsilon}$.  Thus
$Z_{-2\varepsilon}=\Sigma(K_-)\cong\Sigma(K_+)=Z_0$, and Theorem~\ref{thm:monotone} shows that $\tg$ is
null-homologous.  The lift of the crossing arc to $\Sigma(K_-)$ is the dual knot of $M(\beta)$, and by
Lemma~\ref{lem:isotropic}(c) its homology class has the same order as $[\tg]$.  This proves~(a), and~(b) follows
from~(a) and Corollary~\ref{cor:LMnugatory}.
\end{proof}

Theorem~\ref{thm:ccc} contains Theorems~\ref{thm:A} and~\ref{thm:B} for crossing changes, and the following corollary
contains Corollary~\ref{cor:A} for crossing changes; generalized crossing changes are treated in
Section~\ref{subsec:gcc}.

\begin{corollary}\label{cor:ccc}
Every knot whose branched double cover is an L-space satisfies the cosmetic crossing conjecture.  In particular, every
quasi-alternating knot, and hence every alternating knot, satisfies the cosmetic crossing conjecture.
\end{corollary}

\begin{proof}
If a crossing change on a knot $K$ yields a knot $K'$ isotopic to $K$ as an oriented or as an unoriented knot, then
$\Sigma(K')\cong\Sigma(K)$, and Theorem~\ref{thm:ccc}(b) shows that the crossing change is nugatory.  The branched
double covers of quasi-alternating knots are L-spaces
\cite[Proposition~3.3]{OS05}, and alternating knots are quasi-alternating \cite[Lemma~3.2]{OS05}.
\end{proof}

By a theorem of Ito, every prime knot with at most ten crossings satisfies the cosmetic crossing conjecture, with the
possible exception of $\kn{10}{66}$, $\kn{10}{87}$, $\kn{10}{98}$, $\kn{10}{129}$ and $\kn{10}{147}$
\cite[Corollary~1.11]{ItoCW}.
The first three are alternating, and the other two are quasi-alternating \cite{KnotInfo}.  Hence
Corollary~\ref{cor:ccc} gives Corollary~\ref{cor:B}, which we restate; the knots with at most twelve
crossings are discussed in Appendix~\ref{app:census}.

\begin{corollary}\label{cor:ten}
Every prime knot with at most ten crossings satisfies the cosmetic crossing conjecture.
\end{corollary}

Part~(a) of Theorem~\ref{thm:ccc} holds for every knot.  Together with Proposition~\ref{prop:LMnugatory}, it reduces
the cosmetic crossing conjecture for a knot $K$ to a question about null-homologous knots in $Y=\Sigma(K)$: if every
null-homologous knot $\kappa\subset Y$ with $Y_{1/2}(\kappa)\cong Y$ or $Y_{-1/2}(\kappa)\cong Y$ bounds an embedded
disk, then $K$ satisfies the conjecture.  For L-spaces this is Theorem~\ref{thm:gainullin}.

The hypothesis of Theorem~\ref{thm:ccc} concerns orientation-preserving homeomorphisms, and it cannot be weakened to
$\Sigma(K_-)\cong\pm\Sigma(K_+)$.  In Example~\ref{ex:pretzel}, a crossing change turns the knot $6_1$ into its mirror
image, whose branched double cover is $-\Sigma(6_1)$, and the lift of the crossing arc has order three.

\begin{remark}\label{rem:ccclinks}
Theorem~\ref{thm:ccc} extends to links.  Let $L_+$ and $L_-$ be links that differ at a crossing $c$, and suppose that
$\Sigma(L_-)\cong\Sigma(L_+)$ and $\det L_+\neq0$.  By Remark~\ref{rem:links}, the Montesinos trick and
Lemma~\ref{lem:LM} hold for $L_\pm$, and the proof of Theorem~\ref{thm:ccc}(a) shows that the lift of the crossing arc
is null-homologous.  Suppose moreover that $\Sigma(L_+)$ is an L-space and that the two strands at $c$ belong to the
same component.  Then the crossing change is nugatory.  Indeed, the proof of \cite[Proposition~3.3]{LM} applies the equivariant
Dehn lemma and Smith theory to the covering involution, whose fixed point set is the lift of the link, and otherwise
uses only the slopes of Lemma~\ref{lem:LM} and the fact that, of the two curves on the boundary of a ball $B$ around
$c$ whose lifts have distance one from both $\alpha$ and $\beta$, only the crossing circle is null-homologous in the
complement of the link.  This fact holds here.  Indeed, the crossing circle bounds a disk in $\partial B$ that meets the link
in a point where it enters $B$ and a point where it leaves $B$, whereas the other curve bounds a disk in $\partial B$
that meets the link in the two points where it enters $B$.  If $m$ is a meridian of the component that contains the
two strands at $c$, the two curves are therefore homologous to $0$ and to $\pm2m\neq0$ in the complement of the
link.  In particular, links with a connected alternating diagram admit no
cosmetic crossing change between two strands of the same component.
\end{remark}

\begin{remark}\label{rem:ccczhs}
As in \cite[Theorem~3.5]{LM}, Theorem~\ref{thm:ccc} and Corollary~\ref{cor:ccc} hold for knots $K$ in an integral
homology sphere $X$ such that the double cover $\Sigma(X,K)$ of $X$ branched along $K$ is an L-space, where a crossing
change is defined as $\pm1$-surgery along a crossing circle.  Indeed, the Montesinos trick applies in any integral homology
sphere, $\Sigma(X,K)$ is a rational homology sphere, and the proof of \cite[Proposition~3.3]{LM} applies to any closed
$3$-manifold with an involution \cite[proof of Theorem~3.5]{LM}.
\end{remark}

\subsection{Generalized crossing changes}\label{subsec:gcc}
By \cite[Remark~2.7]{Rog}, the generalized crossing change of order $|q|$ given by $1/q$-surgery on a crossing circle
$c$ is the two-strand $(-2q)$-twist determined by $c$, and every two-strand twist of nonzero even order arises in this
way.  We use the notation of Section~\ref{subsec:montesinos}: $\tg\subset\Sigma(K)$ is the lift of a crossing arc
$\gamma$, with exterior $M$ and meridian $\tilde\mu$.

\begin{theorem}\label{thm:gcc}
Let $K'$ be obtained from a knot $K$ by a generalized crossing change, and suppose that $\Sigma(K')\cong\Sigma(K)$; for
instance, $K'$ is isotopic to $K$, as an oriented or as an unoriented knot.
\begin{enumerate}[label=(\alph*)]
\item The lift $\tg$ of the crossing arc to $\Sigma(K)$ is null-homologous.
\item If $\Sigma(K)$ is an L-space, then the generalized crossing change is nugatory; in particular, $K'$ is isotopic
to $K$.
\end{enumerate}
\end{theorem}

\begin{proof}
The generalized crossing change is the two-strand $n$-twist determined by the crossing circle for a nonzero even
integer $n$, so $K'=K_n$ in the notation of~\eqref{eq:twistfamily}.  Since $|H_1(\Sigma(K'))|=|H_1(\Sigma(K))|$, we
have $\det K'=\det K$, and Lemma~\ref{lem:gcc} shows that $\Delta(\tilde\mu,\lambda_M)=1$ and
$\Sigma(K')=M(\tilde\mu+\varepsilon n\lambda_M)$.  By Lemma~\ref{lem:isotropic}(a), the knot $\tg$ is isotropic.  In the
notation~\eqref{eq:family}, with $\kappa=\tg$, $\mu=\tilde\mu$ and $\lambda=\lambda_M$, we have
$\Sigma(K')=M(\mu_{-\varepsilon n})=Z_{-\varepsilon n}$.  Thus $Z_{-\varepsilon n}\cong Z_0$, and since $n\neq0$,
Theorem~\ref{thm:monotone} shows that $\tg$ is null-homologous.  This proves~(a), and~(b) follows from~(a) and
Corollary~\ref{cor:rogersgcc}.
\end{proof}

Theorem~\ref{thm:gcc} proves Theorems~\ref{thm:A} and~\ref{thm:B}, and the following corollary, together with
Corollary~\ref{cor:ccc}, proves Corollary~\ref{cor:A}.

\begin{corollary}\label{cor:gcc}
Every knot whose branched double cover is an L-space satisfies the generalized cosmetic crossing conjecture
\cite[Conjecture~1.1]{Rog}: a generalized crossing change on it that yields an isotopic knot is nugatory.  In
particular, this holds for every quasi-alternating knot, and hence for every alternating knot.
\end{corollary}

\begin{proof}
If a generalized crossing change on a knot $K$ yields a knot $K'$ isotopic to $K$, then $\Sigma(K')\cong\Sigma(K)$, and
Theorem~\ref{thm:gcc}(b) applies.  The second assertion follows as in the proof of Corollary~\ref{cor:ccc}.
\end{proof}

\subsection{Band surgeries and mutants}\label{subsec:bands}
The same argument applies to non-coherent band surgeries, which by \cite[Remark~4.2]{Rog} are the two-strand twists of
order one of Section~\ref{subsec:montesinos}.

\begin{theorem}\label{thm:band}
Let $K'$ be obtained from a knot $K\subset S^3$ by a two-strand twist of order one, for instance by a non-coherent band
surgery, and let $\tg\subset\Sigma(K)$ be the lift of a twisting arc in its twisting disk.  Suppose that $\Sigma(K')\cong\Sigma(K)$; for
instance, $K'$ is isotopic to $K$.
\begin{enumerate}[label=(\alph*)]
\item The knot $\tg$ is null-homologous.
\item If $\Sigma(K)$ is an L-space, then the twist is weakly nugatory.
\end{enumerate}
\end{theorem}

\begin{proof}
By~\eqref{eq:rogers}, $\Sigma(K')=M(\tilde\mu')$ with $\Delta(\tilde\mu,\tilde\mu')=1$, where $M$ is the exterior and
$\tilde\mu$ the meridian of $\tg$.  Since $\det K'=\det K$, Lemma~\ref{lem:band} gives $\Delta(\tilde\mu,\lambda_M)=1$
and $\tilde\mu'=\pm(\tilde\mu+\delta\lambda_M)$.  Hence $\tg$ is isotropic, and in the notation~\eqref{eq:family},
with $\kappa=\tg$, $\mu=\tilde\mu$ and $\lambda=\lambda_M$, we have $\Sigma(K')=Z_{-\delta}\cong Z_0$.  By
Theorem~\ref{thm:monotone}, $\tg$ is null-homologous, and~(b) follows from Corollary~\ref{cor:rogers}.
\end{proof}

Following Ichihara and Jong \cite{IJ}, a band surgery on a knot $K$ is \emph{purely cosmetic} if it yields a knot
isotopic to $K$ as an unoriented knot, and \emph{chirally cosmetic} if it yields the mirror image of $K$.  Both are non-coherent,
since a coherent band surgery on a knot yields a link with two components.  Every knot
admits purely cosmetic band surgeries that are not nugatory \cite[Example~4.1 and the discussion following
it]{Rog}, and Rogers asked whether there
is a nontrivial knot in $S^3$ with a purely cosmetic band surgery that is not weakly nugatory
\cite[Question~4.7]{Rog}.  Theorem~\ref{thm:band} answers this question in the negative for knots with L-space
branched double covers, and together with the following corollary it proves Theorem~\ref{thm:D}.

\begin{corollary}\label{cor:band}
If $\Sigma(K)$ is an L-space, then every purely cosmetic band surgery on $K$ is weakly nugatory.
\end{corollary}

Again the orientation-preserving hypothesis cannot be dropped.  The manifold $X$ of Example~\ref{ex:BHW} is the
exterior of the lift of the twisting arc of a chirally cosmetic band surgery on the two-bridge knot $9_{27}$
\cite[Section~5.1]{IJ}, and this lift is not null-homologous.

Two knots are \emph{mutants} if they are related by a sequence of Conway mutations.  Viro showed that mutants have
homeomorphic branched double covers \cite[Theorem~1]{Vir76}; the homeomorphism constructed there is
orientation-preserving, since the gluing maps are lifts of rotations that are isotopic to the identity.  For
alternating knots the converse holds by a theorem of Greene \cite[Theorem~1.1]{Gre13}: two reduced alternating diagrams
represent knots $K$ and $K'$ with $\Sigma(K')\cong\Sigma(K)$ if and only if they are related by mutation.

\begin{corollary}\label{cor:mutants}
Let $K$ and $K'$ be knots with $\Sigma(K')\cong\Sigma(K)$, for instance mutants, and suppose that $\Sigma(K)$ is an
L-space.  If $K'$ is obtained from $K$ by a crossing change, then $K'$ is isotopic to $K$.  In particular, two mutant
knots with L-space branched double covers that are not isotopic have Gordian distance at least two.
\end{corollary}

\begin{proof}
This follows from Theorem~\ref{thm:ccc}(b).
\end{proof}

\begin{remark}\label{rem:bandmutants}
Corollary~\ref{cor:mutants} does not extend to non-coherent band surgeries, and in Theorem~\ref{thm:band}(b) the knot
$K'$ need not be isotopic to $K$, even if the twist is nugatory.  Let $J$ be the knot $8_{17}$, which is alternating
and not invertible \cite{KnotInfo}, and let $rJ$ be $J$ with its orientation reversed.  Let $S$ be a sphere that meets
$K=J\#J$ in two points and separates its two summands, and let $D\subset S$ be a disk that contains these two points.
Then $c=\partial D$ is a twisting circle for $K$ with twisting disk $D$, and the twists along $D$ are nugatory, since
$c$ also bounds the disk $S\setminus\operatorname{int}D$, which is disjoint from $K$.  A two-strand twist of odd order
along $D$ amounts to cutting $S^3$ along $S$ and regluing by a rotation of $S$ through $\pi$ that exchanges the two
points of $K\cap S$, so it yields $K'=J\#rJ$.  The oriented prime factors of $K$, $rK$ and $K'$ are
$\{J,J\}$, $\{rJ,rJ\}$ and $\{J,rJ\}$, so by the uniqueness of prime decompositions of oriented knots, $K'$ is not
isotopic to $K$, even as an unoriented knot.  On the other hand, $\Sigma(K')\cong\Sigma(J)\#\Sigma(rJ)=\Sigma(K)$,
which is an L-space, since $K$ is alternating.  Both knots have reduced alternating diagrams, so by the theorem of
Greene quoted above they are mutants.

The same circle $c$ also bounds a twisting disk along which the twists of odd order do not change $K$.  Indeed, let
$D_0$ be a disk that meets $K$ in two points and lies on the boundary of a small ball meeting $K$ in an unknotted arc.
A twist of odd order along $D_0$ replaces this arc by another unknotted arc with the same endpoints, so it yields $K$.
The circles $c$ and $\partial D_0$ bound disks in $S^3\setminus K$, so an isotopy of $S^3$ fixing $K$ carries
$\partial D_0$ to $c$, and $D_0$ to a twisting disk bounded by $c$.  Hence, for odd $n$, the two-strand $n$-twist
depends on the twisting disk and not only on its boundary, in contrast to \cite[Remark~2.7]{Rog}.
\end{remark}

\subsection{Knots determined by their complements}\label{subsec:complement}
Let $\kappa$ be a knot in a closed $3$-manifold $Y$, with exterior $M$ and meridian $\mu$, and recall the notion of
a knot determined by its complement and Conjecture~\ref{conj:kc} from Section~\ref{subsec:introkc}.  If
$M(\eta)\not\cong Y$ for every slope $\eta\neq\mu$, then $\kappa$ is determined by its complement.  Indeed, a homeomorphism
between the complements of $\kappa'$ and $\kappa$ induces one between the exteriors \cite[Section~1]{Gai}, which
carries the meridian of $\kappa'$ to a slope $\eta$ on $\partial M$ with $M(\eta)\cong Y$, so $\eta=\mu$ and the
homeomorphism extends over $Y$.  Theorem~\ref{thm:monotone} gives the following, which is Theorem~\ref{thm:E}(a).

\begin{theorem}\label{thm:complement}
Let $\kappa$ be a knot in a rational homology sphere $Y$ whose homology class is isotropic and nonzero.  Then
$M(\eta)\not\cong Y$ for every slope $\eta\neq\mu$.  In particular, $\kappa$ is determined by its complement.
\end{theorem}

\begin{proof}
By Lemma~\ref{lem:isotropic}(a), the rational longitude $\lambda_M$ is a framing of $\kappa$, so
Theorem~\ref{thm:monotone} applies.
\end{proof}

The theorem does not extend to orientation-reversing homeomorphisms.  In Example~\ref{ex:BHW}, a nontrivial surgery
on a knot of order seven in $L(49,18)$ yields $-L(49,18)$.  For lens spaces, Theorem~\ref{thm:complement} completes
the work of Gainullin; the following corollary contains Theorem~\ref{thm:E}(b) for lens spaces, and
Theorem~\ref{thm:elliptic} below contains the general case.

\begin{corollary}\label{cor:lens}
Let $L=L(p,q)$ be a lens space, $p\ge1$, and let $K\subset L$ be a knot that does not bound a disk and whose exterior
is not a solid torus.  Then no
nontrivial surgery on $K$ yields $L$ by an orientation-preserving homeomorphism.  In particular, in a lens
space, every knot whose exterior is irreducible and not a solid torus is determined by its complement.
\end{corollary}

\begin{proof}
Let $M$ be the exterior and $\mu$ the meridian of $K$, and suppose that $M(\eta)\cong L$ for a slope $\eta\neq\mu$.

If $K$ is null-homologous, write $\eta$ as $a/b$ with respect to the Seifert framing.  Then
$|H_1(M(\eta))|=|a|\,|H_1(L)|$, so $|a|=1$, and $L_{1/b}(K)\cong L\cong L_{1/b}(U)$ for an unknot $U\subset L$ (with
$a=1$ after changing the signs of $a$ and $b$).  Since lens spaces are L-spaces, Theorem~\ref{thm:gainullin} shows that
$K$ is isotopic to $U$, which is excluded.

Suppose that $K$ is not null-homologous.  If its class is isotropic, then Theorem~\ref{thm:complement} contradicts
$M(\eta)\cong L$, so we may assume that it is not.  The exterior $M$ is irreducible.  Indeed, since $L$ is irreducible, an
essential sphere in $M$ would bound a ball in $L$ containing $K$, and $K$ would be null-homologous.

Suppose first that $M$ is Seifert fibered.  Filling $M$ along the fiber slope gives a manifold whose fundamental group
is the orbifold fundamental group of the base of $M$, which is infinite since $M$ is not a solid torus.  Hence $\mu$
is not the fiber slope, the Seifert fibration of $M$ extends over $L$, and $K$ is a fiber.  Since $\pi_1(L)$ is
finite and cyclic, the base of this fibration of $L$ is the projective plane with at most one cone point or the
sphere with at most two cone points.  In the first case, let $\alpha$ be the order of the cone point, or $\alpha=1$ if
there is none.  The
standard presentation of $\pi_1(L)$ shows that $|H_1(L)|=4\alpha$ and that the class of a regular fiber, which is
conjugate to its inverse, has order two.  The element of order two of the cyclic group $H_1(L)$ is isotropic, since
$4$ divides $|H_1(L)|$, so $K$ is not a regular fiber.  Hence $K$ is the exceptional fiber and $M$ is Seifert fibered
over the M\"obius band without exceptional fibers; that is, $M$ is the twisted $I$-bundle over the Klein bottle,
which is also Seifert fibered over the disk with two exceptional fibers of order two.  Applying the above to this
fibration, we may assume that the base of the fibration of $L$ is the sphere.  Since $M$ is not a solid torus, it
has at least two exceptional fibers, so $L$ has exactly two and $K$ is a regular fiber; that is, $K$ is a nontrivial
cable of the core of a Heegaard solid torus of $L$.  We show that $K$ is then isotropic, contrary to our assumption.

In the notation of Geiges and Lange \cite[Section~2]{GeL18}, with respect to the orientation of $L$, the exterior $M$
is Seifert fibered over the disk with two exceptional fibers of types $(\alpha_1,\beta_1)$ and $(\alpha_2,\beta_2)$,
where $\alpha_1,\alpha_2\ge2$.  Let $h$ be the fiber slope on $\partial M$ and $s$ the boundary slope of a section.
The slopes at distance one from $h$ are the slopes $s+kh$ with $k\in\Z$, and the filling along $s+kh$ is
\[
M_k=M\bigl(0;(\alpha_1,\beta_1),(\alpha_2,\beta_2),(1,k)\bigr)=M\bigl(0;(\alpha_1,\beta_1),(\alpha_2,\beta_2+k\alpha_2)\bigr).
\]
Since $M(\eta)\cong L$, the argument above, applied to $\eta$ in place of $\mu$, shows that $\eta$ is not the fiber
slope and that the core of the filling along $\eta$ is a regular fiber.  Hence $\mu=s+nh$ and $\eta=s+n'h$ for
integers $n\neq n'$.  By \cite[Theorem~4.4]{GeL18}, $M_k=L(p_k,q_k)$ for every integer $k$, where
\[
p_k=\alpha_1\beta_2+\alpha_2\beta_1+k\alpha_1\alpha_2,\qquad q_k=\alpha_1(\beta'_2+k\alpha'_2)+\beta_1\alpha'_2,
\]
the integers $\alpha'_2$ and $\beta'_2$ satisfy $\alpha_2\beta'_2-\beta_2\alpha'_2=1$, and $L(p,q)$ stands for $L(-p,-q)$
if $p<0$.  In particular, $q_k\alpha_2=\alpha_1+\alpha'_2p_k$.

Since $|p_n|=|H_1(L)|=|p_{n'}|$ and $p_{n'}-p_n=(n'-n)\alpha_1\alpha_2\neq0$, we have $p_{n'}=-p_n$.  Write
$g=\gcd(\alpha_1,\alpha_2)$, $\alpha_i=ga_i$ and $p_n=gp'$.  Then $2p'=(n-n')ga_1a_2$, and
$p'=a_1\beta_2+a_2\beta_1+nga_1a_2$ is prime to $a_1a_2$, since $\gcd(a_1,a_2)=\gcd(a_1,\beta_1)=\gcd(a_2,\beta_2)=1$;
hence $a_1a_2$ divides $2$.  In particular, one of $a_1$ and $a_2$ is $1$, so $g=\min(\alpha_1,\alpha_2)\ge2$.
Dividing $q_k\alpha_2=\alpha_1+\alpha'_2p_k$ by $g$ for $k=n$ and $k=n'$ gives
$q_na_2\equiv a_1\equiv q_{n'}a_2\pmod{p'}$, so $q_{n'}\equiv q_n\pmod{p'}$.  Now
$M_{n'}=L(-p_n,q_{n'})=L(p_n,-q_{n'})$, and by the classification of lens spaces \cite[Section~3.2]{GeL18}, an
orientation-preserving homeomorphism between $M_{n'}$ and $M_n$ exists only if $-q_{n'}\equiv q_n$ or
$-q_{n'}q_n\equiv1\pmod{p_n}$.  Reducing modulo $p'$, we obtain $2q_n\equiv0$ or $q_n^2\equiv-1\pmod{p'}$.  In the
first case $p'$ divides $2$, since $q_n$ is prime to $p_n$.  In the second case $a_1^2+a_2^2\equiv0\pmod{p'}$, so $p'$
divides $2$ if $a_1=a_2=1$, and $p'$ divides $5$ if $\{a_1,a_2\}=\{1,2\}$.

Finally, write $\lambda_M=\pm(xs+yh)$ with $\gcd(x,y)=1$.  By~\eqref{eq:watson},
$|p_n+(k-n)\alpha_1\alpha_2|=|p_k|=c_M\,|y-kx|$ for every integer $k$ with $s+kh\neq\lambda_M$, so
$c_M=\gcd(p_n,\alpha_1\alpha_2)=g\gcd(g,p')$, and $\Delta(\mu,\lambda_M)=|p_n|/c_M=|p'|/\gcd(g,p')$.  If $|p'|=1$,
this is $1$.  If $|p'|=2$, then $a_1a_2$ is odd, so $a_1=a_2=1$, and $|n-n'|\,g=4$ with $g\ge2$ shows that $g$ is even.  If
$|p'|=5$, then $\{a_1,a_2\}=\{1,2\}$, and $|n-n'|\,g=5$ with $g\ge2$ shows that $g=5$.  In each case
$\Delta(\mu,\lambda_M)=1$, so $K$ is isotropic by Lemma~\ref{lem:isotropic}(a).

If $M$ is not Seifert fibered, then, since $\pi_1(M(\mu))$ and $\pi_1(M(\eta))$ are cyclic, the cyclic surgery theorem
of Culler, Gordon, Luecke and Shalen \cite{CGLS} gives $\Delta(\mu,\eta)\le1$, so $\Delta(\mu,\eta)=1$.  Since
$|H_1(M(\eta))|=|H_1(M(\mu))|$, it follows from~\eqref{eq:watson} that $\Delta(\eta,\lambda_M)=\Delta(\mu,\lambda_M)$;
writing $\mu=m\rho+r\lambda_M$ and $\eta=\pm(m\rho+s\lambda_M)$ with $\rho\cdot\lambda_M=1$ as in the proof of
Lemma~\ref{lem:band}, we obtain $1=m\,|r-s|$, so $\Delta(\mu,\lambda_M)=m=1$.  Thus $K$ is isotropic by
Lemma~\ref{lem:isotropic}(a), which we have excluded.
\end{proof}

Corollary~\ref{cor:lens} was known for non-hyperbolic knots \cite[Corollary~1.2]{Mat10}, when $|H_1(L)|$ is
square-free \cite[Corollary~8.3]{Gai}, when $|H_1(L)|=d^2p'$ with $p'$ square-free and $d\in\{2,3,6\}$ prime to
$p'$ \cite[Theorem~1.7(ii)]{ItoCW}, and for knots that represent a generator of $H_1(L)$ \cite[Theorem~1.2]{IS22}.  In the proof of \cite[Corollary~8.3]{Gai},
Gainullin
showed that if a knot in $L(p,q)$ that is not null-homologous and has irreducible exterior that is not Seifert fibered
admits a nontrivial surgery yielding $L(p,q)$, then the surgery slope is integral and the winding number $w$ satisfies
$p\mid w^2$; this expresses the isotropy of $K$ used above.  The orientation-preserving hypothesis cannot be
dropped.  Matignon showed that a hyperbolic knot in a lens space with a nontrivial surgery that yields the lens space
and preserves a Heegaard solid torus is not determined by its complement \cite[Theorem~1.1]{Mat17}, and he constructed
infinitely many non-hyperbolic knots in lens spaces with such surgeries \cite[Theorems~1.2 and~6.1]{Mat17}.  In
\cite{Mat17}, a \emph{meridional cosmetic slope} of a knot $K$ in a lens space $L$ is a slope $r\neq\mu_K$ on the
exterior $M_K$ of $K$ such that $M_K(r)$ is homeomorphic to $L$, with no condition on orientation
\cite[Section~2.1, p.~395]{Mat17}, whereas the homeomorphism between knot complements in the definition of a knot
determined by its complement is assumed to preserve orientation \cite[Section~1, p.~391]{Mat17}.  What is proved is
that no self-homeomorphism of $M_K$ carries $\mu_K$ to $r$ \cite[Section~6.2, p.~428]{Mat17}.  To compare the two
notions, let $\phi\colon M_K(r)\to L$ be a homeomorphism, and let $K'\subset L$ be the image of the core of the
attached solid torus; then $\phi$ restricts to a homeomorphism from $M_K$ onto the exterior of $K'$ that carries $r$
to the meridian of $K'$.  By Corollary~\ref{cor:lens}, $\phi$ reverses orientation, in particular for the hyperbolic
knots of \cite[Theorem~1.1]{Mat17} and in Example~\ref{ex:BHW}; for the non-hyperbolic knots, Matignon observed that
there is no orientation-preserving homeomorphism between $M_K(r)$ and $L$ \cite[p.~1900]{Mat10}.  For the hyperbolic
knots, every homeomorphism $g$ from the exterior of $K'$ onto $M_K$ reverses orientation as well: otherwise $g$ would
carry the meridian of $K'$ to a slope $s$ with $M_K(s)\cong L$, so $s=\mu_K$ by Corollary~\ref{cor:lens}, and
$g\circ\phi|_{M_K}$ would be a self-homeomorphism of $M_K$ that carries $r$ to $\mu_K$, which
\cite[Section~6.2]{Mat17} excludes.  The same argument shows that no homeomorphism of $L$ carries $K'$ to $K$.  Hence
\cite[Theorem~1.1]{Mat17} holds in the unoriented sense: the complements of $K$ and $K'$ are homeomorphic, but only by
orientation-reversing homeomorphisms.  In the oriented sense of Section~\ref{subsec:introkc}, these knots are
determined by their complements.  The case of the projective plane in the proof of Corollary~\ref{cor:lens} gives further
examples.  By \cite[Proposition~3.1]{GeL21}, $L(8,3)=M(-1;(2,1))$ in the notation of \cite[Section~2]{GeL18}, and a
surgery on a regular fiber of this Seifert fibration, which is isotropic, yields
$M(-1;(2,1),(1,-1))=M(-1;(2,-1))=L(8,5)\cong-L(8,3)$.

For the other elliptic manifolds, knots that are not isotropic need, besides results from the literature, the
following lemma on two fillings at distance three.

\begin{lemma}\label{lem:distance3}
Let $M$ be a compact, connected $3$-manifold with torus boundary such that $H_1(M;\Q)\cong\Q$, and let $\mu$ and
$\eta$ be slopes on $\partial M$ with $\Delta(\mu,\eta)=\Delta(\mu,\lambda_M)=3$ and $|H_1(M(\mu))|=|H_1(M(\eta))|$.  Then
the linking forms of $M(\mu)$ and $M(\eta)$ are not isomorphic.  In particular, $M(\eta)\not\cong M(\mu)$.
\end{lemma}

\begin{proof}
Choose a slope $\nu$ with $\Delta(\mu,\nu)=1$, and orient $\mu$ and $\nu$ so that $\mu\cdot\nu=1$.  Then
$\eta=\pm(a\mu+3\nu)$ for an integer $a$ prime to $3$, and for an integer $b$ with $3b-a=\pm1$ the slope
$\xi=b\mu+\nu$ satisfies $\Delta(\xi,\mu)=\Delta(\xi,\eta)=1$.  Since $\Delta(\xi,\mu)=1\neq\Delta(\lambda_M,\mu)$, we have
$\xi\neq\lambda_M$, so $Y=M(\xi)$ is a rational homology sphere by~\eqref{eq:watson}, and $\mu$ and $\eta$ are framings of
its dual knot $\kappa$.  Present $Y=\partial X_N$ as in Section~\ref{subsec:isotropic}, with $\kappa\subset S^3\setminus
L$ and linking vector $v$.  Then $M(\mu)=\partial X_B$ and $M(\eta)=\partial X_{B'}$ with
\[
B=\begin{pmatrix}N&v\\ v^T&c\end{pmatrix},\qquad B'=\begin{pmatrix}N&v\\ v^T&c'\end{pmatrix},
\]
where $c$ and $c'$ are the coefficients of $\mu$ and $\eta$ with respect to the Seifert framing of $\kappa$ in $S^3$.
Thus $|c-c'|=\Delta(\mu,\eta)=3$, and after exchanging $\mu$ and $\eta$ we may assume that $c'=c+3$.  Then $\det
B'=\det B+3\det N$, and since $|\det B'|=|\det B|$ and $\det N\neq0$, we have $\det B'=-\det B$.

For a nonsingular symmetric integral matrix $A$, the linking form of $\partial X_A$ is the form $(x,y)\mapsto-x^TA^{-1}y$
on $\coker A$ \cite[Section~3]{GL78}.  Since $3$ is odd, $B$ is congruent over the ring $\Z_3$ of $3$-adic integers to a
diagonal matrix.  Indeed, an entry of minimal $3$-adic valuation can be moved to the diagonal, replacing a basis
vector $e_i$ by $e_i+e_j$ if necessary, and the rank-one summand spanned by the corresponding basis vector then splits
off orthogonally.  Collecting the diagonal entries of valuation $k$, we obtain
$B\cong\bigoplus_{k\ge0}3^kU_k$ with $U_k$ unimodular over $\Z_3$ of rank $r_k$.  Hence the $3$-primary part of the
linking form of $\partial X_B$ is isomorphic to the orthogonal sum over $k\ge1$ of the forms $-3^{-k}U_k^{-1}$ on
$(\Z/3^k)^{r_k}$.  By a theorem of Wall \cite[Lemma~(8) and Theorem~(4)]{Wal63}, its isomorphism class determines,
and is determined by, the ranks $r_k$ and the Legendre symbols $\epsilon_k=\bigl(\frac{(-1)^{r_k}\det U_k}{3}\bigr)$ for
$k\ge1$; here, for a unit $w$ of $\Z_3$, $(\frac w3)$ denotes the Legendre symbol of its residue modulo $3$.  Write
$\det B=3^tu$ with $u$ prime to $3$.  Then the product of the determinants of the $U_k$ is congruent to
$u$ modulo squares of units of $\Z_3$, so $\prod_{k\ge0}\bigl(\frac{\det U_k}{3}\bigr)=\bigl(\frac u3\bigr)$.

Define $U'_k$, $r'_k$, $\epsilon'_k$ and $u'$ for $B'$ in the same way.  Since $B'\equiv B\pmod3$, the reductions of
$U_0$ and $U'_0$ modulo $3$ are both isomorphic to the nondegenerate quotient of the same form over the field with
three elements, so $r_0=r'_0$ and $\bigl(\frac{\det U'_0}{3}\bigr)=\bigl(\frac{\det U_0}{3}\bigr)$.  Since $\det
B'=-\det B$, we have $u'=-u$, and since $-1$ is not a square modulo $3$, it follows that
\[
\prod_{k\ge1}\Bigl(\frac{\det U'_k}{3}\Bigr)=-\prod_{k\ge1}\Bigl(\frac{\det U_k}{3}\Bigr).
\]
Hence $r_k\neq r'_k$ or $\epsilon_k\neq\epsilon'_k$ for some $k\ge1$, and in both cases the $3$-primary parts of the
linking forms of $M(\mu)$ and $M(\eta)$ are not isomorphic.  Since an orientation-preserving homeomorphism induces an
isomorphism of linking forms, $M(\eta)\not\cong M(\mu)$.
\end{proof}

The hypothesis $\Delta(\mu,\lambda_M)=3$ cannot be dropped.  For instance, if $M$ is the exterior of a knot that
bounds a disk in a rational homology sphere $Y$, then $\Delta(\mu,\lambda_M)=1$, and $M(\mu)$ and $M(\mu+3\lambda_M)$ are both
homeomorphic to $Y$ by orientation-preserving homeomorphisms.

\begin{theorem}\label{thm:elliptic}
Let $Y$ be an elliptic $3$-manifold, and let $K\subset Y$ be a knot that does not bound a disk and whose exterior is not
a solid torus.  Then no nontrivial surgery on $K$ yields $Y$ by an orientation-preserving homeomorphism.  In
particular, in $Y$, every knot whose exterior is irreducible and not a solid torus is determined by its complement.
\end{theorem}

\begin{proof}
If $Y$ is a lens space, this is Corollary~\ref{cor:lens}, so we assume that it is not.  We use the following properties
of $Y$.  It is irreducible, since its universal cover is $S^3$ \cite[Proposition~1.6]{Hat}.  It is an L-space
\cite[Proposition~2.3]{OSlens}.  Its fundamental group is finite, so every torus in $Y$ is compressible.  Finally, $Y$
admits no orientation-reversing homeomorphism \cite[Proposition~1.1]{McC02}.

Let $M$ be the exterior and $\mu$ the meridian of $K$, let $\lambda=\lambda_M$, and suppose that $M(\eta)\cong Y$ for a
slope $\eta\neq\mu$.  As in the proof of Corollary~\ref{cor:lens}, Theorem~\ref{thm:gainullin} and
Theorem~\ref{thm:complement} show that $K$ is not isotropic, so $m=\Delta(\mu,\lambda)\ge2$ by
Lemma~\ref{lem:isotropic}(a), and $M$ is irreducible.  Since $|H_1(M(\eta))|=|H_1(M(\mu))|$, \eqref{eq:watson} gives
$\Delta(\eta,\lambda)=m$.  Writing $\mu=m\rho+r\lambda$ and $\eta=\pm(m\rho+s\lambda)$ with $\rho\cdot\lambda=1$, we obtain
$d=\Delta(\mu,\eta)=m\,|r-s|$, and if $m=2$, then $r$ and $s$ are odd, so that $4\mid d$.  Thus $m$ divides $d$, and
$d\ge3$.  Since $M$ is irreducible and not a solid torus, its boundary is incompressible, so by Thurston's
hyperbolization theorem \cite[Theorem~2.3]{Thu82}, in the form stated in \cite[Section~1.0]{CGLS}, $M$ is Seifert
fibered, contains an essential torus, or is hyperbolic.  We derive
a contradiction in each case, using the following two observations.

\emph{(i) Let $M_0\subset Y$ be the exterior of a knot in $Y$.  Suppose that $M_0$ is Seifert fibered and not a solid
torus, and that $M_0(\mu_0)\cong M_0(\eta_0)\cong Y$ for slopes $\mu_0$ and $\eta_0$.  Then $\mu_0=\eta_0$.}  Filling
$M_0$ along the fiber slope gives a manifold whose fundamental group is infinite, as in the proof of
Corollary~\ref{cor:lens}, so neither $\mu_0$ nor $\eta_0$ is the fiber slope, and a Seifert fibration of $M_0$ extends to
Seifert fibrations $\mathcal F_{\mu_0}$ of $M_0(\mu_0)$ and $\mathcal F_{\eta_0}$ of $M_0(\eta_0)$, whose base surfaces
are both obtained from that of $M_0$ by capping off its boundary circle.  Let $f\colon M_0(\eta_0)\to M_0(\mu_0)$ be an
orientation-preserving diffeomorphism.  The Seifert fibrations $f(\mathcal F_{\eta_0})$ and $\mathcal F_{\mu_0}$ of
$M_0(\mu_0)$ have the same base surface, so they are isomorphic by \cite[Theorem~2.3]{Hat}, which is proved for
fibrations over $S^2$ with three multiple fibers in \cite{Orl72}.  Among the exceptions listed in \cite[Theorem~2.3]{Hat}, the only
manifolds with finite fundamental group that are not lens spaces are those of case~(d), whose two
Seifert fibrations have the base surfaces $S^2$ and $\mathbb{RP}^2$.  An isomorphism between $f(\mathcal F_{\eta_0})$ and
$\mathcal F_{\mu_0}$ is a diffeomorphism of $M_0(\mu_0)\cong Y$, so it preserves orientation.  Hence there is an
orientation-preserving diffeomorphism $M_0(\eta_0)\to M_0(\mu_0)$ that carries $\mathcal F_{\eta_0}$ to $\mathcal
F_{\mu_0}$.  In the notation of \cite[Section~2.1]{Hat}, write $M_0=M(\pm g,1;\alpha_1/\beta_1,\dots,\alpha_k/\beta_k)$
with respect to a section; then $M_0(\mu_0)$ and $M_0(\eta_0)$ are the fibrations
$M(\pm g,0;\alpha_1/\beta_1,\dots,\alpha_k/\beta_k,\alpha/\beta)$, where $\alpha/\beta$ is the slope of $\mu_0$,
respectively $\eta_0$, with respect to the section and the fiber.  By \cite[Proposition~2.1(b)]{Hat}, the two sums
$\sum_i\alpha_i/\beta_i+\alpha/\beta$ agree, so $\mu_0=\eta_0$.

\emph{(ii) If $M_0\subset Y$ is the exterior of a knot in $Y$, $M_0$ is hyperbolic, and $M_0(\mu_0)\cong
M_0(\eta_0)\cong Y$, then $\Delta(\mu_0,\eta_0)\le3$.}  This is the finite filling theorem of Boyer and Zhang \cite[Theorem~1.1]{BZ01}, since
$\pi_1(Y)$ is finite.

If $M$ is Seifert fibered, then~(i) contradicts $\mu\neq\eta$.  If $M$ is hyperbolic, then $d\le3$ by~(ii), so $d=3$
and, since $m\ge2$ divides $d$, also $m=3$; this contradicts Lemma~\ref{lem:distance3}.

Suppose that $M$ contains an essential torus $T$.  It compresses in $M(\mu)$ and in $M(\eta)$, and since $d\ge2$, the
theorem of Culler, Gordon, Luecke and Shalen \cite[Theorem~2.0.1]{CGLS} shows that $T$ and $\partial M$ cobound a cable
space $C\subset M$, of type $(p,q)$ with $q\ge2$.  Let $M'$ be the closure of $M\setminus C$; then $\partial M'=T$ is
incompressible in $M'$, and $M'$ is irreducible.  Since $T$ compresses in $M(\mu)=M'\cup C(\mu)$, it compresses in
$C(\mu)$, and likewise in $C(\eta)$.  By \cite[Lemma~7.2]{Gor83}, each of $\mu$ and $\eta$ is therefore either the
cabling slope of $C$ or a slope at distance one from it, and in the second case the filling of $C$ is a solid torus.
Since $\Delta(\mu,\eta)\ge3$, neither $\mu$ nor $\eta$ is the cabling slope, so $C(\mu)$ and $C(\eta)$ are solid tori.
Hence $M(\mu)=M'(\mu')$ and $M(\eta)=M'(\eta')$ for
slopes $\mu'$ and $\eta'$ on $T$ with $\Delta(\mu',\eta')=q^2d\ge12$ \cite[Section~2.7]{CGLS}, see also
\cite[Corollary~7.3]{Gor83}.
Thus $M'$ is the exterior of the core of $C(\mu)$, a knot in $Y$; it is not a solid torus, and $M'(\mu')\cong
M'(\eta')\cong Y$.  If $M'$ is Seifert fibered or hyperbolic, this contradicts~(i) or~(ii).  Otherwise $M'$ contains an
essential torus, and we repeat the argument with $M'$ in place of $M$, which gives tori $T'$, $T''$, \dots\ in $M$.  These
tori are disjoint, and each of them is incompressible in $M$ by \cite[Section~1.2, property~(5)]{Hat}, since it is incompressible in the
piece that contains it and the tori found before it are incompressible on both sides.
Cutting $M$ along them gives cable spaces, which are not products $T^2\times[0,1]$ since $q\ge2$, and one further piece
with connected boundary.  Hence no two of these tori are parallel in $M$, since an incompressible torus in
$T^2\times[0,1]$ is parallel to a boundary component.  By a finiteness theorem of Haken \cite[Proposition~1.7]{Hat},
the number of these tori is
bounded, so the process stops, and we obtain a contradiction.
\end{proof}

\begin{remark}\label{rem:sister}
Lemma~\ref{lem:distance3} cannot be omitted from the proof of Theorem~\ref{thm:elliptic}, and the orientation-preserving
hypothesis cannot be dropped for elliptic manifolds that are not lens spaces.  Let $M$ be the exterior of the
figure-eight sister knot, the knot $K$ in $L(5,1)$ obtained from the Whitehead link by $5$-surgery on one of its
components, and let $\mu$ and $\lambda$ be the standard meridian and longitude of the other component.  By
\cite[Example~11.7]{BZ01}, $M$ is hyperbolic, and the fillings along $3\mu+\lambda$ and $3\mu+2\lambda$ have finite
fundamental groups of $T$-type, that is, isomorphic to $T(8,3^k)\times\Z/j$ with $k\ge1$ and $\gcd(6,j)=1$, whose first
homology is $\Z/3^k\oplus\Z/j$.  Since the components of the Whitehead link have linking number zero, $K$ is
null-homologous in $L(5,1)$, so $\lambda$ is the rational longitude of $M$ and $H_1(M(p\mu+q\lambda))\cong\Z/5\oplus\Z/p$.
Hence both fillings have first homology $\Z/15$, their fundamental groups are isomorphic to $T(8,3)\times\Z/5$, and the
two slopes have distance three from each other and from $\lambda$.  Elliptic manifolds with isomorphic non-cyclic
fundamental groups are homeomorphic \cite[Section~1]{McC02}, and by Lemma~\ref{lem:distance3} the two fillings are not
homeomorphic by an orientation-preserving homeomorphism.  Therefore $M(3\mu+2\lambda)\cong-M(3\mu+\lambda)$, and the dual
knot $\kappa$ in $Y=M(3\mu+\lambda)$ is a hyperbolic knot with a nontrivial surgery yielding $-Y$, so the first
assertion of Theorem~\ref{thm:elliptic} fails for orientation-reversing homeomorphisms.  Nevertheless, $\kappa$ is
determined by its complement even in the unoriented sense.  SnapPy \cite{SnapPy} finds an orientation-reversing
isometry of $M$ that acts on slopes by $p\mu+q\lambda\mapsto\pm\bigl(p\mu+(p-q)\lambda\bigr)$; it exchanges
$3\mu+\lambda$ and $3\mu+2\lambda$, and it extends to an orientation-reversing homeomorphism from $M(3\mu+\lambda)$ to
$M(3\mu+2\lambda)$, which gives a second proof of the last homeomorphism.  By \cite[Example~11.7]{BZ01}, the finite
fillings of $M$ are those along $\mu$, $\mu+\lambda$, $2\mu+\lambda$, $3\mu+\lambda$ and $3\mu+2\lambda$, and only the
last two have non-cyclic fundamental groups.  Now let $\kappa'\subset Y$ be a knot, and let $h$ be a homeomorphism, of
either orientation, from the exterior of $\kappa'$ onto $M$.  Then $h$ carries the meridian of $\kappa'$ to
$3\mu+\lambda$ or to $3\mu+2\lambda$, and after composing $h$ with the isometry above if necessary, to $3\mu+\lambda$.
Then $h$ extends to a homeomorphism of $Y$ that carries $\kappa'$ to $\kappa$, and this homeomorphism preserves
orientation, since $Y$ admits no orientation-reversing homeomorphism \cite[Proposition~1.1]{McC02}.
\end{remark}

Before Theorem~\ref{thm:elliptic}, the oriented knot complement conjecture was known in every elliptic $3$-manifold for
null-homologous knots \cite[Theorem~8.2]{Gai}, hence for all knots in the Poincar\'e homology sphere, for all knots when
the manifold is obtained by surgery on a knot in $S^3$ and has square-free first homology \cite[Theorem~1.7(iii)]{ItoCW},
and for hyperbolic knots that represent a generator of a cyclic first homology group of order divisible by a prime at
least $7$ \cite[Theorem~1.2(2)]{IS22}.  For non-hyperbolic knots, the conclusion of Theorem~\ref{thm:elliptic} is stated
in \cite[Corollary~1.4]{Mat10}; observation~(i) and the toroidal case above give a proof for elliptic manifolds that does
not depend on \cite{Mat10}.

\subsection*{Use of AI tools}

Claude Opus 5.5 (Anthropic) was used to check formulas numerically and
symbolically. It was also used to carry out the computations in
Appendix~\ref{app:census}. The author verified all the results and takes full
responsibility for the content of this paper.

\subsection*{Acknowledgments}
This work is supported by the Institute for Basic Science (IBS) under Project No.\ IBS-R003-D1. The author would like to thank National Center for Theoretical Sciences for hospitality during an important stage of this work. 
\appendix

\section{A lattice-theoretic proof for alternating knots}\label{app:alternating}

For alternating knots, part~(a) of Theorem~\ref{thm:ccc}, and with it Corollary~\ref{cor:ccc}, can also be deduced
from Donaldson's diagonalization theorem in place of Theorem~\ref{thm:monotone}.  The input from an alternating diagram
is a pair of negative definite $2$-handlebodies, bounded by the branched double cover with its two orientations, whose
intersection lattices are the negatives of the cut lattices of the two Tait graphs.  We use the notation of
Section~\ref{subsec:isotropic}, and we identify a symmetric integral $n\times n$ matrix with the lattice $\Z^n$
carrying the corresponding form; for instance, a vector of square $-1$ of $N$ is a vector $y\in\Z^n$ with
$y^TNy=-1$.

\subsection{Chains of \texorpdfstring{$2$}{2}-handles}\label{subsec:chains}
Let $Y=\partial X_N$ and let $\kappa\subset S^3\setminus L$ be an isotropic knot in $Y$ with linking vector $v$, so
that $x_0=v^TN^{-1}v$ is an integer by Lemma~\ref{lem:presentation}.  For an integer $q\ge1$, let
\[
C_q=\begin{pmatrix}-1&1&&\\1&-2&\ddots&\\&\ddots&\ddots&1\\&&1&-2\end{pmatrix},
\qquad
B_q(v)=\begin{pmatrix}N&v\,u^T\\ u\,v^T&C_q+x_0\,u\,u^T\end{pmatrix},
\]
where $C_q$ is a $q\times q$ matrix and $u=(1,0,\dots,0)^T\in\Z^q$.  Then $B_q(v)$ is the linking matrix of the framed
link consisting of $L$, of $\kappa$ with framing $x_0-1$, and of a chain of $q-1$ unknots with framing $-2$, the first
a meridian of $\kappa$ and each of the others a meridian of the previous one.

\begin{lemma}\label{lem:chain}
For every $q\ge1$, the manifold $Y_{x_0-1/q}(\kappa)$ obtained by surgery on $\kappa$ with coefficient $x_0-1/q$ is
the boundary of $X_{B_q(v)}$, and $\det B_q(v)=(-1)^q\det N$.  If $N$ is negative definite, then so is $B_q(v)$.
\end{lemma}

\begin{proof}
The chain of $q-1$ unknots with framing $-2$ has continued fraction
\[
-2-\cfrac{1}{-2-\cfrac{1}{\ddots-\cfrac{1}{-2}}}=-\frac{q}{q-1}\qquad(q\ge2),
\]
so slam-dunks \cite{GS} turn the framing $x_0-1$ of $\kappa$ into the coefficient $x_0-1+(q-1)/q=x_0-1/q$.  This
proves the first assertion.  For the others, let $e_1,\dots,e_n,g,h_1,\dots,h_{q-1}$ be the basis of
$H_2(X_{B_q(v)})$ given by the $2$-handles, with $g$ attached along $\kappa$.  Put $\hat v=N^{-1}v\in\Q^n$ and
$g'=g-\sum_{i=1}^n\hat v_ie_i\in H_2(X_{B_q(v)};\Q)$.  Then $g'\cdot e_i=0$ for all $i$, $g'\cdot g'=-1$ as in the proof
of Lemma~\ref{lem:W}(c), and $g'\cdot h_1=1$ when $q\ge2$.  Hence, over $\Q$, the form $B_q(v)$ is the orthogonal sum of $N$ and of
the form $C_q$ on the basis $g',h_1,\dots,h_{q-1}$.  The vectors
\[
f_j=g'+h_1+\dots+h_{j-1}\qquad(1\le j\le q)
\]
satisfy $f_j\cdot f_j=-1$ and $f_i\cdot f_j=0$ for $i\neq j$, and they form a basis of the span of
$g',h_1,\dots,h_{q-1}$ over $\Z$.  Thus $C_q$ is congruent to $-I_q$.  Since the change of basis from
$e_1,\dots,e_n,g,h_1,\dots,h_{q-1}$ to $e_1,\dots,e_n,g',h_1,\dots,h_{q-1}$ has determinant one, $\det B_q(v)=\det
N\cdot\det C_q=(-1)^q\det N$, and $B_q(v)$ is negative definite whenever $N$ is.
\end{proof}

With respect to the framing $\lambda=\lambda_M$, which has coefficient $x_0$, the slope with coefficient $x_0-1/q$ is
$\mu-q\lambda$; thus $\partial X_{B_q(v)}$ is the manifold $Z_q$ of~\eqref{eq:family}.

\begin{lemma}\label{lem:minusone}
Let $N$ be negative definite with no vector of square $-1$, let $v$ be an integral vector with $x_0\in\Z$, and let
$q\ge1$.  If $[v]\neq0$ in $\coker N$, then $B_q(v)$ has no vector of square $-1$.
\end{lemma}

\begin{proof}
We use the notation of the proof of Lemma~\ref{lem:chain}.  An element of $H_2(X_{B_q(v)})$ has the form
$\omega=\sum_iy_ie_i+a\,g+\sum_jb_jh_j$ with $y\in\Z^n$ and $a,b_1,\dots,b_{q-1}\in\Z$.  Substituting
$g=g'+\sum_i\hat v_ie_i$, we obtain
\[
\omega\cdot\omega=(y+a\,\hat v)^TN\,(y+a\,\hat v)+Q(a,b),
\]
where $Q(a,b)$ denotes the value of the form $C_q$ on the vector with entries $a,b_1,\dots,b_{q-1}$.  Both terms are
nonpositive, and $Q(a,b)$ is an integer.  Suppose that $\omega\cdot\omega=-1$.  If $(a,b)=0$, then $y^TNy=-1$, which
is excluded.  Otherwise $Q(a,b)\le-1$, so $Q(a,b)=-1$ and $y+a\,\hat v=0$.  The vectors of square $-1$ of $C_q\cong-I_q$
are the vectors $\pm f_j$, and each of them has coefficient $\pm1$ on $g'$.  Hence $a=\pm1$ and $\hat v=\mp y\in\Z^n$,
so
$v=N\hat v$ and $[v]=0$.
\end{proof}

\subsection{Goeritz lattices}\label{subsec:goeritz}
Let $G$ be a finite connected graph, possibly with multiple edges, and orient its edges.  The \emph{star} of a vertex
$u$ is the vector
\[
s_u=\sum_{e\colon u\to\,\cdot}e-\sum_{e\colon\cdot\,\to u}e\in\Z^{E(G)},
\]
where the first sum runs over the edges leaving $u$ and the second over the edges entering it.  The \emph{cut lattice}
$C(G)\subset\Z^{E(G)}$ is the sublattice spanned by the stars, with the inner product of $\Z^{E(G)}$; its isometry type
does not depend on the orientations of the edges.  The \emph{flow lattice} $F(G)$ is the orthogonal complement of
$C(G)$ in $\Z^{E(G)}$.  For a lattice $\Lambda$ we write $-\Lambda$ for the same group with the negated form.  We shall
use the following elementary facts.
\begin{enumerate}[label=(\roman*)]
\item The stars satisfy $\sum_us_u=0$.  If $G$ has no loops, then $s_u\cdot s_u=\deg u$ for every vertex $u$, and
hence
\[
\sum_{u\in V(G)}s_u\cdot s_u=2\,|E(G)|.
\]
\item Every element of $C(G)$ is the coboundary $\delta a$ of a function $a\colon V(G)\to\Z$, where $(\delta
a)_e=a(u)-a(u')$ for an edge $e$ from $u$ to $u'$.  Consequently $C(G)$ contains a vector of square $1$ if and only if
$G$ has a bridge, since $\delta a$ has a single nonzero entry, equal to $\pm1$, exactly when $a$ is constant on the two
components of the complement of a bridge and its values on them differ by one.
\item If $G$ is a plane graph and $G^*$ is its planar dual, then $F(G)=C(G^*)$ under the identification
$E(G)=E(G^*)$, since, suitably oriented, the boundaries of the faces of $G$ are the stars of $G^*$, and they span the lattice
of integral flows of $G$.
\end{enumerate}

Let $\mathcal D$ be a reduced alternating diagram of a knot $K$ with $n$ crossings, and fix a checkerboard coloring of
its complementary regions.  The two Tait graphs of $\mathcal D$ have the regions of one color as vertices and the
crossings as edges; they are planar duals of each other, and each of them has $n$ edges.  Since $\mathcal D$ is
reduced, no circle in the plane meets $\mathcal D$ in a single crossing and nowhere else, so neither graph has a loop
or a bridge.  The two graphs have $n+2$ vertices in total, one for each region.  Since all crossings of an alternating
diagram have the same incidence number with respect to the coloring, the Goeritz matrix of $\mathcal D$ formed on the
regions of one color
\cite[\S1]{GL78} is, up to sign, the Gram matrix of the cut lattice of the Tait graph on these regions, with respect to
the stars of all but one vertex; for this reason we refer to the cut lattices of the Tait graphs as Goeritz lattices.

\begin{lemma}[Ozsv\'ath--Szab\'o]\label{lem:goeritz}
Let $Y=\Sigma(K)$.  The Tait graphs of $\mathcal D$ can be labeled $G_+$ and $G_-$ so that $Y$ and $-Y$ bound negative
definite $2$-handlebodies $X_+$ and $X_-$ with intersection lattices
\[
H_2(X_+)\cong-C(G_+),\qquad H_2(X_-)\cong-C(G_-).
\]
\end{lemma}

\begin{proof}
By \cite[Lemma~3.6]{OS05} and its proof, $\Sigma(K)$ bounds a negative definite $4$-manifold, obtained from $B^4$ by
attaching $2$-handles after canceling handle pairs, whose intersection form is the form $Q$ of \cite{OS05}, that is,
minus the flow lattice of the black graph of $\mathcal D$ for the coloring convention of \cite[Figure~3]{OS05}.
By~(iii), this is minus the cut lattice of the other Tait graph, which we call $G_+$.  Applying the same result to the
mirror image of $\mathcal D$, whose black graph is $G_+$, we find that $\Sigma(\overline K)=-Y$ bounds a negative
definite $2$-handlebody with intersection lattice $-F(G_+)=-C(G_-)$.
\end{proof}

By~(ii), neither $H_2(X_+)$ nor $H_2(X_-)$ contains a vector of square $-1$.  Since the cut lattice of a connected
graph has rank one less than its number of vertices, $\rk H_2(X_+)+\rk H_2(X_-)=n$.  Apart from~(ii), the only
property of cut lattices that we need is the following count.

\begin{lemma}\label{lem:counting}
Let $G$ be a connected graph without loops, and let $\phi\colon C(G)\to\Z^m$ be an isometric embedding into the
standard lattice.  Then at least $m-|E(G)|$ of the standard basis vectors of $\Z^m$ are orthogonal to $\phi(C(G))$.
\end{lemma}

\begin{proof}
Let $\xi_1,\dots,\xi_m$ be the standard basis, and write $\phi(s_u)=\sum_i\phi_{ui}\,\xi_i$.  By~(i),
$\sum_u\phi_{ui}=0$ for every $i$.  Hence, if the column $(\phi_{ui})_u$ is nonzero, it has a positive and a negative
entry, and $\sum_u\phi_{ui}^2\ge2$.  Again by~(i),
\[
\sum_{i=1}^m\sum_{u\in V(G)}\phi_{ui}^2=\sum_{u\in V(G)}\phi(s_u)\cdot\phi(s_u)=2\,|E(G)|,
\]
so at most $|E(G)|$ columns are nonzero.  If the column of $i$ is zero, then $\xi_i$ is orthogonal to every
$\phi(s_u)$, and hence to $\phi(C(G))$.
\end{proof}

\subsection{The lattice argument}\label{subsec:latticeproof}

\begin{theorem}\label{thm:lattice}
Let $K$ be an alternating knot and $Y=\Sigma(K)$.  Let $\kappa\subset Y$ be an isotropic knot with exterior $M$ and
meridian $\mu$, and orient $\mu$ and $\lambda=\lambda_M$ so that $\mu\cdot\lambda=1$ on $\partial\nu(\kappa)$.  If
$M(\mu+q\lambda)\cong Y$ for some nonzero integer $q$, then $\kappa$ is null-homologous.
\end{theorem}

\begin{proof}
If $K$ is trivial, then $Y=S^3$ and there is nothing to prove, so we may assume that $K$ is nontrivial.  Replacing
$K$ by its mirror image, which is again alternating and has branched double cover $-Y$, reverses the orientations of
$Y$ and of $M$; to keep $\mu\cdot\lambda=1$ we then replace $\lambda$ by $-\lambda$, which replaces $q$ by $-q$.  We may
therefore assume that $q<0$, so that $M(\mu+q\lambda)=M(\mu-|q|\lambda)$ is the manifold $Z_{|q|}$
of~\eqref{eq:family}.  Fix a reduced alternating diagram $\mathcal D$ of $K$ with $n\ge1$ crossings, and let $X_Y$ and
$X_{-Y}$ be the $2$-handlebodies of Lemma~\ref{lem:goeritz} with $\partial X_Y=Y$ and $\partial X_{-Y}=-Y$.  Let $N$ be
the intersection matrix of $X_Y$ in the basis given by its $2$-handles, and let $v$ be a linking vector of $\kappa$
with respect to the corresponding framed link.  By Lemma~\ref{lem:presentation}, $x_0=v^TN^{-1}v$ is an integer, and by
Lemma~\ref{lem:chain} and the remark following it, $\partial X'=Z_{|q|}\cong Y$, where $X'=X_{B_{|q|}(v)}$ is negative
definite.

Choose an orientation-preserving diffeomorphism $\partial X'\to Y$, which exists because every homeomorphism between
$3$-manifolds is isotopic to a diffeomorphism, and glue $X'$ to $X_{-Y}$ along it.  The result is a closed, oriented,
smooth $4$-manifold $\widehat X$.  Since $X'$ and $X_{-Y}$ are simply connected, so is $\widehat X$, and
$H_2(\widehat X)$ is free.  Since $Y$ is a rational homology sphere, there is an orthogonal decomposition
\[
H_2(\widehat X;\Q)=H_2(X';\Q)\oplus H_2(X_{-Y};\Q),
\]
so $\widehat X$ is negative definite, and
\[
b_2(\widehat X)=\bigl(\rk H_2(X_Y)+|q|\bigr)+\rk H_2(X_{-Y})=n+|q|.
\]
By Donaldson's theorem \cite[Theorem~1]{Don}, $H_2(\widehat X)$ has a basis $\xi_1,\dots,\xi_{n+|q|}$ with
$\xi_i\cdot\xi_j=-\delta_{ij}$.  The inclusion of $X_{-Y}$ induces an isometric embedding of $H_2(X_{-Y})\cong-C(G)$,
where $G$ is a Tait graph of $\mathcal D$, into $H_2(\widehat X)$.  Applying Lemma~\ref{lem:counting} to the negatives
of these forms, we find that at least $|q|\ge1$ of the vectors $\xi_i$ are orthogonal to $H_2(X_{-Y})$.

We claim that every such $\xi_i$ lies in $H_2(X')$.  By excision and Poincar\'e--Lefschetz duality,
\begin{gather*}
H_3(\widehat X,X')\cong H_3(X_{-Y},\partial X_{-Y})\cong H^1(X_{-Y})=0,\\
H_2(\widehat X,X')\cong H_2(X_{-Y},\partial X_{-Y})\cong H^2(X_{-Y}),
\end{gather*}
and $H^2(X_{-Y})$ is free.  By the exact sequence of the pair $(\widehat X,X')$, the map $H_2(X')\to H_2(\widehat X)$
is therefore injective with torsion-free cokernel.  Hence $H_2(X')$ is a primitive sublattice of $H_2(\widehat X)$,
and by the orthogonal decomposition above it is the orthogonal complement of $H_2(X_{-Y})$.  This proves the claim.
Consequently $H_2(X')$, whose intersection matrix is $B_{|q|}(v)$, contains a vector of square $-1$.  Since $N$ has no
vector of square $-1$, Lemma~\ref{lem:minusone} shows that $[v]=0$; that is, $\kappa$ is null-homologous in $Y$.
\end{proof}

If $\Sigma(K_-)\cong\Sigma(K_+)$ for an alternating knot $K_+$, then Lemma~\ref{lem:LM} shows that the lift $\tg$ of
the crossing arc is isotropic and that $M(\beta)=M(\alpha+2\varepsilon\lambda_M)\cong M(\alpha)$, so
Theorem~\ref{thm:lattice}, applied with $\mu=\alpha$ and $q=2\varepsilon$, gives Theorem~\ref{thm:ccc}(a) for
alternating knots.  Together with Corollary~\ref{cor:LMnugatory}, this proves Corollary~\ref{cor:ccc} for alternating
knots without Theorem~\ref{thm:monotone}.  Suppose now that a generalized crossing change given by a nonzero even
number $n$ of half-twists turns an alternating knot $K$ into a knot $K_n$ with $\Sigma(K_n)\cong\Sigma(K)$.  Then
$\det K_n=\det K$, so Lemma~\ref{lem:gcc} identifies the twisting slope with the rational longitude, the lift of the
twisting arc is isotropic by Lemma~\ref{lem:isotropic}(a), and $\Sigma(K_n)=M(\tilde\mu+\varepsilon n\lambda_M)$.
Theorem~\ref{thm:lattice}, applied with $q=\varepsilon n$, shows that this lift is null-homologous, and
Corollary~\ref{cor:rogersgcc} shows that the twist is nugatory.  This proves Corollary~\ref{cor:A} for alternating
knots without Theorem~\ref{thm:monotone}.  The argument needs negative definite $2$-handlebodies bounded by
$\Sigma(K)$ and by $-\Sigma(K)$ whose intersection lattices have total rank $n$, such that one lattice has no vector of
square $-1$ and the other is the negative of the cut lattice of a graph with $n$ edges.  For alternating knots
these are the Goeritz lattices of the two checkerboard surfaces of a reduced alternating diagram, one positive and one
negative definite, and by a theorem of Greene a knot bounds both a positive definite and a negative definite spanning
surface only if it is alternating \cite[Theorem~1.1]{Gre17}.

\section{Knots with at most twelve crossings}\label{app:census}

In November 2026, \KI{} \cite{KnotInfo} listed the cosmetic crossing conjecture as open for $214$ prime knots with at most twelve crossings.
Of these, $68$ are alternating and satisfy both the ordinary and the generalized conjecture by Corollary~\ref{cor:A}.
For the other $146$ we compared $\det K$ with the rank of $\HFhat(\Sigma(K);\Z/2)$ computed from the plat presentations in Zhan's bordered Floer program \cite{LOT,Zha}.
The ranks agree with the distributed data (\nolinkurl{data/output_12.txt}, commit \texttt{e96f345}).
Equality holds for the $62$ knots in Table~\ref{tab:resolved}, so Corollary~\ref{cor:gcc} applies over $\Z/2$ as explained after Theorem~\ref{thm:gainullin}.
For the $84$ knots in Table~\ref{tab:open} the rank exceeds the determinant.
Thus both conjectures hold for $130$ of the $214$ knots. For each of the remaining $84$, the crossing arc of any cosmetic generalized crossing change would have null-homologous lift by Theorem~\ref{thm:gcc}(a).

Among these $130$, the Smith normal form of $V+V^T$, using the Seifert matrix $V$ in \KI{}, gives square-free order in each invariant factor of $H_1(\Sigma(K))$ in exactly $14$ cases \cite[Theorem~1.2]{LM}.
Five of the other $116$ have irreducible Alexander polynomial, so Boninger's criterion also excludes cosmetic generalized crossing changes for these knots \cite[Theorem~1.2]{Bon24}.

There is a second way to verify $57$ of the $62$ entries in Table~\ref{tab:resolved}.
For $43$, the knot is quasi-alternating or its reduced Khovanov homology over $\Z/2$ has rank $\det K$ \cite{KnotInfo,OS05}.
For the $14$ Montesinos knots in Table~\ref{tab:seifert}, the table gives the normalized Seifert invariants of the cover up to orientation \cite[Theorem~A.6]{BS}.
By \cite[Theorems~1.1 and~2.1]{LS}, it suffices to exclude coprime integers $m>a>0$ with $s_1<a/m<1-s_2$ and $s_3<1/m$.
Every row has $s_3\ge1/5$, leaving only $m=2,3,4$ to check. None gives such a fraction.
The remaining five are $\kn{11n}{126}$, $\kn{12n}{260}$, $\kn{12n}{265}$, $\kn{12n}{442}$ and $\kn{12n}{810}$.

\begin{table}[ht]
\caption{The 62 non-alternating knots with at most twelve crossings, listed as open in \KI{}, whose branched double
covers are L-spaces with $\Z/2$ coefficients.}\label{tab:resolved}
{\footnotesize
\setlength{\columnsep}{4pt}
\setlength{\multicolsep}{3pt}
\par\vspace{3pt}\hrule height\heavyrulewidth
\begin{multicols}{10}
\raggedright\setlength{\parindent}{0pt}\setlength{\parskip}{0pt}
$\kn{10}{129}$\par
$\kn{10}{147}$\par
$\kn{11n}{20}$\par
$\kn{11n}{50}$\par
$\kn{11n}{68}$\par
$\kn{11n}{69}$\par
$\kn{11n}{72}$\par
$\kn{11n}{79}$\par
$\kn{11n}{83}$\par
$\kn{11n}{122}$\par
$\kn{11n}{126}$\par
$\kn{11n}{132}$\par
$\kn{11n}{138}$\par
$\kn{12n}{48}$\par
$\kn{12n}{49}$\par
$\kn{12n}{50}$\par
$\kn{12n}{62}$\par
$\kn{12n}{65}$\par
$\kn{12n}{66}$\par
$\kn{12n}{153}$\par
$\kn{12n}{170}$\par
$\kn{12n}{176}$\par
$\kn{12n}{200}$\par
$\kn{12n}{224}$\par
$\kn{12n}{243}$\par
$\kn{12n}{244}$\par
$\kn{12n}{251}$\par
$\kn{12n}{260}$\par
$\kn{12n}{265}$\par
$\kn{12n}{286}$\par
$\kn{12n}{310}$\par
$\kn{12n}{322}$\par
$\kn{12n}{331}$\par
$\kn{12n}{360}$\par
$\kn{12n}{364}$\par
$\kn{12n}{365}$\par
$\kn{12n}{380}$\par
$\kn{12n}{393}$\par
$\kn{12n}{412}$\par
$\kn{12n}{413}$\par
$\kn{12n}{421}$\par
$\kn{12n}{422}$\par
$\kn{12n}{440}$\par
$\kn{12n}{442}$\par
$\kn{12n}{444}$\par
$\kn{12n}{502}$\par
$\kn{12n}{503}$\par
$\kn{12n}{523}$\par
$\kn{12n}{553}$\par
$\kn{12n}{556}$\par
$\kn{12n}{597}$\par
$\kn{12n}{637}$\par
$\kn{12n}{667}$\par
$\kn{12n}{702}$\par
$\kn{12n}{726}$\par
$\kn{12n}{731}$\par
$\kn{12n}{782}$\par
$\kn{12n}{799}$\par
$\kn{12n}{810}$\par
$\kn{12n}{813}$\par
$\kn{12n}{841}$\par
$\kn{12n}{876}$\par
\end{multicols}
\hrule height\heavyrulewidth\par}
\end{table}

\begin{table}[ht]
\caption{The 84 non-alternating knots in the November 2026 \KI{} list for which $\operatorname{rk}\HFhat(\Sigma(K);\Z/2)>\det K$. The present criterion leaves these knots undecided.}\label{tab:open}
{\footnotesize
\setlength{\columnsep}{4pt}
\setlength{\multicolsep}{3pt}
\par\vspace{3pt}\hrule height\heavyrulewidth
\begin{multicols}{10}
\raggedright\setlength{\parindent}{0pt}\setlength{\parskip}{0pt}
$\kn{11n}{34}$\par
$\kn{11n}{39}$\par
$\kn{11n}{42}$\par
$\kn{11n}{45}$\par
$\kn{11n}{67}$\par
$\kn{11n}{73}$\par
$\kn{11n}{97}$\par
$\kn{11n}{151}$\par
$\kn{11n}{152}$\par
$\kn{12n}{23}$\par
$\kn{12n}{28}$\par
$\kn{12n}{31}$\par
$\kn{12n}{34}$\par
$\kn{12n}{51}$\par
$\kn{12n}{56}$\par
$\kn{12n}{57}$\par
$\kn{12n}{63}$\par
$\kn{12n}{82}$\par
$\kn{12n}{87}$\par
$\kn{12n}{96}$\par
$\kn{12n}{110}$\par
$\kn{12n}{118}$\par
$\kn{12n}{124}$\par
$\kn{12n}{129}$\par
$\kn{12n}{132}$\par
$\kn{12n}{149}$\par
$\kn{12n}{195}$\par
$\kn{12n}{203}$\par
$\kn{12n}{210}$\par
$\kn{12n}{214}$\par
$\kn{12n}{217}$\par
$\kn{12n}{221}$\par
$\kn{12n}{225}$\par
$\kn{12n}{231}$\par
$\kn{12n}{232}$\par
$\kn{12n}{256}$\par
$\kn{12n}{257}$\par
$\kn{12n}{258}$\par
$\kn{12n}{264}$\par
$\kn{12n}{267}$\par
$\kn{12n}{268}$\par
$\kn{12n}{279}$\par
$\kn{12n}{292}$\par
$\kn{12n}{293}$\par
$\kn{12n}{313}$\par
$\kn{12n}{321}$\par
$\kn{12n}{332}$\par
$\kn{12n}{374}$\par
$\kn{12n}{382}$\par
$\kn{12n}{386}$\par
$\kn{12n}{394}$\par
$\kn{12n}{404}$\par
$\kn{12n}{411}$\par
$\kn{12n}{430}$\par
$\kn{12n}{432}$\par
$\kn{12n}{457}$\par
$\kn{12n}{464}$\par
$\kn{12n}{473}$\par
$\kn{12n}{483}$\par
$\kn{12n}{519}$\par
$\kn{12n}{535}$\par
$\kn{12n}{554}$\par
$\kn{12n}{575}$\par
$\kn{12n}{594}$\par
$\kn{12n}{638}$\par
$\kn{12n}{644}$\par
$\kn{12n}{650}$\par
$\kn{12n}{655}$\par
$\kn{12n}{665}$\par
$\kn{12n}{676}$\par
$\kn{12n}{680}$\par
$\kn{12n}{689}$\par
$\kn{12n}{691}$\par
$\kn{12n}{692}$\par
$\kn{12n}{750}$\par
$\kn{12n}{764}$\par
$\kn{12n}{801}$\par
$\kn{12n}{808}$\par
$\kn{12n}{809}$\par
$\kn{12n}{812}$\par
$\kn{12n}{815}$\par
$\kn{12n}{830}$\par
$\kn{12n}{851}$\par
$\kn{12n}{868}$\par
\end{multicols}
\hrule height\heavyrulewidth\par}
\end{table}

\begin{table}[ht]
\caption{The fourteen Montesinos knots in Table~\ref{tab:resolved} to which \cite{OS05} does not apply,
with their rational tangles $(r_1,r_2,r_3)$ in the notation of \KI{}, where $r_i=\beta_i/\alpha_i$ in lowest terms
with $\alpha_i\ge2$ and $\det K=\alpha_1\alpha_2\alpha_3|r_1+r_2+r_3|$, and the numbers $(s_1,s_2,s_3)$ with
$\Sigma(K)\cong\pm M(-1;s_1,s_2,s_3)$.}\label{tab:seifert}
{\footnotesize
\begin{tabular}{lll@{\qquad}lll}
\toprule
$K$ & $(r_1,r_2,r_3)$ & $(s_1,s_2,s_3)$ & $K$ & $(r_1,r_2,r_3)$ & $(s_1,s_2,s_3)$\\
\midrule
$\kn{11n}{20}$ & $(1/2, -3/5, 3/7)$ & $(1/2, 3/7, 2/5)$ & $\kn{12n}{244}$ & $(1/2, -2/3, 7/17)$ & $(1/2, 7/17, 1/3)$\\
$\kn{11n}{79}$ & $(1/2, -2/3, 4/9)$ & $(1/2, 4/9, 1/3)$ & $\kn{12n}{251}$ & $(1/2, -2/3, 5/11)$ & $(1/2, 5/11, 1/3)$\\
$\kn{11n}{138}$ & $(2/3, -1/3, -4/7)$ & $(4/7, 1/3, 1/3)$ & $\kn{12n}{310}$ & $(1/2, -5/7, 3/7)$ & $(1/2, 3/7, 2/7)$\\
$\kn{12n}{65}$ & $(1/2, -3/5, 4/9)$ & $(1/2, 4/9, 2/5)$ & $\kn{12n}{322}$ & $(2/3, -3/5, -2/7)$ & $(3/5, 1/3, 2/7)$\\
$\kn{12n}{153}$ & $(1/2, -2/3, 5/17)$ & $(1/2, 1/3, 5/17)$ & $\kn{12n}{502}$ & $(1/2, -4/7, 1/5)$ & $(1/2, 3/7, 1/5)$\\
$\kn{12n}{200}$ & $(1/2, -2/3, 4/15)$ & $(1/2, 1/3, 4/15)$ & $\kn{12n}{503}$ & $(1/2, -4/7, 3/7)$ & $(1/2, 3/7, 3/7)$\\
$\kn{12n}{243}$ & $(1/2, -2/3, 3/13)$ & $(1/2, 1/3, 3/13)$ & $\kn{12n}{523}$ & $(2/3, -1/4, -4/7)$ & $(4/7, 1/3, 1/4)$\\
\bottomrule
\end{tabular}}
\end{table}

\FloatBarrier

\end{document}